\documentclass[a4paper,11pt]{amsart}%

\usepackage{amsmath}
\usepackage{amsfonts}
\usepackage{amssymb}
\usepackage{hyperref}
\usepackage{enumitem}
\usepackage{xcolor}
\usepackage{graphicx}
\usepackage{epstopdf}
\usepackage{algorithm}
\usepackage{todonotes}

\usepackage{tikz}
\usepackage{tikz-qtree}

\usetikzlibrary{positioning, calc}

\usepackage[noend]{algorithmic}%
\newtheorem{theorem}{Theorem}

\newtheorem{definition}[theorem]{Definition}
\newtheorem{lemma}[theorem]{Lemma}
\theoremstyle{remark}
\newtheorem{example}[theorem]{Example}
\newtheorem{remark}[theorem]{Remark}

\newcommand{\Q}{\mathbb{Q}}
\newcommand{\CC}{\mathbb{C}}
\newcommand{\PP}{\mathbb{P}}
\def\good{/\!\!/}

\begin{document}
\title[Challenges in Developing Open Source Computer Algebra Systems]{Current Challenges in Developing Open Source Computer Algebra Systems}

\author[J.~B\"ohm]{Janko~B\"ohm}
\address{Department of Mathematics\\
University of Kaiserslautern\\
Erwin-Schr\"odinger-Str.\\
67663 Kaiserslautern\\
Germany}
\email{boehm@mathematik.uni-kl.de}

\author[W.~Decker]{Wolfram Decker}
\address{Department of Mathematics\\
University of Kaiserslautern\\
Erwin-Schr\"odinger-Str.\\
67663 Kaiserslautern\\
Germany}
\email{decker@mathematik.uni-kl.de}

\author[S.~Keicher]{Simon~Keicher}
\address{
Departamento de Matematica\\
Facultad de Ciencias Fisicas y Matematicas\\
Universidad de Concepci{\'o}n\\ Casilla 160-C, Concepci{\'o}n, Chile 
}
\email{simonkeicher@googlemail.com}

\author[Y.~Ren]{Yue~Ren}
\address{Department of Mathematics\\
University of Kaiserslautern\\
Erwin-Schr\"odinger-Str.\\
67663 Kaiserslautern\\
Germany}
\email{ren@mathematik.uni-kl.de}

\thanks{The second author acknowledges support from the DFG projects DE 410/8-1 and -2, 
\mbox{DE 410/9-1} \mbox{and -2}, and
from the OpenDreamKit Horizon 2020 European Research Infrastructures project 
($\#$676541).
The third author was supported partially by the DFG project HA 3094/8-1
and by proyecto FONDECYT postdoctorado no 3160016.}

\maketitle

\begin{abstract}
This note is based on the plenary talk given by the second author at MACIS
2015, the Sixth International Conference on Mathematical Aspects of Computer 
and Information Sciences. Motivated by some of the work done within the Priority Programme
SPP 1489 of the German Research Council DFG, we discuss a number of current challenges 
in the development of Open Source computer algebra systems. The main
focus is on algebraic geometry and the system {\sc{Singular}}.
\end{abstract}
\section{Introduction}

The goal of the nationwide Priority Programme SPP 1489 of the German Research Council DFG
is to considerably further the algorithmic and experimental methods in algebraic geometry, 
number theory, and group theory, to combine the different methods where needed, and to 
apply them to central questions in theory and practice. In particular, the programme is meant 
to support the further development of Open Source computer algebra systems which are 
(co-)based in Germany, and which in the framework of different projects may require 
crosslinking on different levels. The cornerstones of the latter are the well-established systems
\textsc{GAP}~\cite{GAP4} (group and representation theory), \textsc{polymake}~\cite{polymake} (polyhedral geometry),
and \textsc{Singular}~\cite{singular} (algebraic geometry, singularity theory, commutative and
non-commutative algebra), together with the newly evolving system \textsc{ANTIC}~\cite{antic}
(number theory), but there are many more systems, libraries, and packages involved (see Section 
\ref{sec systems} for some examples). 

In this note, having the main focus on \textsc{Singular}, we report on some of the 
challenges which we see in this context. These range from reconsidering the efficiency 
of the basic algorithms through parallelization and making abstract concepts constructive 
to facilitating the access to Open Source computer algebra systems. In illustrating the 
challenges, which are discussed in Section \ref{sect:challenges}, we take examples 
from algebraic geometry. In Sections \ref{sec normal} and Section \ref{sec git},
two of the  examples are highlighted in more detail. These are the parallelization 
of the classical Grauert-Remmert type algorithms for normalization 
and the computation of GIT-fans. The latter is a show-case application 
of bringing \textsc{Singular}, \textsc{polymake}, and \textsc{GAP} together.

\section{Seven Challenges}\label{sect:challenges}

\subsection{Reconsidering the Efficiency of the Basic Algorithms}

Motivated by an increasing number of success stories in applying algorithmic 
and experimental methods to algebraic geometry (and other areas of mathematics), 
research projects in this direction become more and more ambitious. This
applies both to the theoretical  level of abstraction and to the practical complexity. 
On the computer algebra side, this not only requires innovative ideas to design 
high-level algorithms, but also to revise the basic algorithms on which the 
high-level algorithms are built. The latter concerns efficiency and applicability.

\begin{example}[The \textsc{Nemo} project] \textsc{Nemo} is a new computer algebra 
package written in the  \textsc{Julia}\footnote{See \url{http://julialang.org}} 
programming language which, in particular, aims
at highly efficient implementations of basic arithmetic and algorithms for number theory
and is connected to the \textsc{ANTIC} project. 
See \url{http://nemocas.org/index.html} for some benchmarks.
\end{example}

In computational algebraic geometry, aside from polynomial factorization, the basic work 
horse is Buchberger's algorithm for computing Gr\"obner bases \cite{zbMATH06070338-ds} 
and, as remarked by Schreyer \cite{SchreyerDiplom-ds} and others, syzygies. While Gr\"obner 
bases are specific sets of generators for ideals and modules which are well-suited for 
computational purposes, the name syzygies refers to the relations on a 
given set of generators. Syzygies carry important geometric information
(see \cite{EisenbudSyz}) and are crucial ingredients in many basic and high-level 
algorithms. Taking syzygies on the syzygies and so forth,
we arrive at what is called a free resolution. Here is a particular simple
example.

\begin{example}[The Koszul complex of Three Variables] In the \textsc{Singular}
session below, we first construct the polynomial ring $R=\mathbb{Q}[x,y,z]$, 
endowed with the degree reverse lexicographical order \texttt{dp}. Then we
compute the successive syzygies on the variables $x,y,z$.

\hspace{-3mm}{ \texttt{$\color{blue}>$ \color{red}ring \color{black}R = 0, (x,y,z), dp;}}

\hspace{-3mm}{ \texttt{$\color{blue}>$ \color{red}ideal \color{black}I = x,y,z;}}

\hspace{-3mm}{ \texttt{$\color{blue}>$ \color{red}resolution \color{black}FI = nres(I,0);}}

\hspace{-3mm}{ \texttt{$\color{blue}>$ \color{red}print\color{black}(FI[2]);}}

\hspace{-3mm}{ \color{blue}\texttt{ \hspace{3mm}  0,-y,-z,}}

\hspace{-3mm}{ \color{blue}\texttt{ \hspace{1mm} -z, x, 0,}}

\hspace{-3mm}{ \color{blue}\texttt{ \hspace{3mm}  y, 0, x}}

\hspace{-3mm}{ \texttt{$\color{blue}>$ \color{red}print\color{black}(FI[3]);}}

\hspace{-3mm}{ \color{blue}\texttt{ \hspace{3mm}  x,}}

\hspace{-3mm}{ \color{blue}\texttt{ \hspace{3mm}  z,}}

\hspace{-3mm}{ \color{blue}\texttt{ \hspace{1mm} -y}}

\end{example}

In the following example, we show how Gr\"obner basis and syzygy
computations fit together to build a more advanced algorithm.

\begin{example}[Parametrizing Rational Curves]\label{ex:para}
We  study a degree-5 curve $C$ in the projective plane which is
visualized as the red curve in Figures \ref{fig adjoints generators} 
and \ref{fig adjoint}.   To begin with, after constructing the polynomial 
ring $R=\mathbb{Q}[x,y,z]$, we enter the homogeneous degree-5 
polynomial $f\in \mathbb{Q}[x,y,z]$ which defines~$C$:

\hspace{-3mm}{ \texttt{$\color{blue}>$ \color{red}ring \color{black}R = 0, (x,y,z), dp;}}

\hspace{-3mm}{ \texttt{$\color{blue}>$ \color{red}poly \color{black}f = x5+10x4y+20x3y2+130x2y3-20xy4+20y5-2x4z-40x3yz }}

\hspace{-3mm}{ \texttt{$\color{blue}\hspace{18mm}$ -150x2y2z-90xy3z-40y4z+x3z2+30x2yz2+110xy2z2+20y3z2;}}

\noindent
Our goal is to check whether  $C$ is rational, and if so, to compute a
rational parametrization. For the first task, recall that an algebraic curve is rational 
if and only if its geometric genus is zero. In the example here,  
this can be easily read off from the genus formula for 
plane curves, taking into account that the degree-5 curve has three 
ordinary double points and one ordinary triple point (see the 
aforementioned visualization). An algorithm for computing the genus
in general, together with 
an algorithm for computing rational parametrizations, is implemented in the
\textsc{Singular} library \texttt{paraplanecurves.lib}~\cite{paraplane}:

\hspace{-3mm}{ \texttt{$\color{blue}>$ \color{red}LIB \color{black} "paraplanecurves.lib";}}

\hspace{-3mm}{ \texttt{$\color{blue}>$ \color{red}genus\color{black}(f);}}

\hspace{-3mm}{ \texttt{\color{blue} \hspace{0.3cm} 0}}

\hspace{-3mm}{ \texttt{$\color{blue}>$ \color{red}paraPlaneCurve\color{black}(f);}}

\noindent 
Rather than displaying the result, we will now show the key steps of the algorithm at work.
The first step is to compute the ideal generated by the adjoint curves of $C$ which,
roughly speaking, are curves which pass with sufficiently high multiplicity
through the singular points of $C$. 
The algorithm for computing the adjoint ideal (see~\cite{adjointideal}) builds on 
algorithms for computing normalization (see Section \ref{sec normal}) or,
equivalently, integral bases (see~\cite{integralbases}). In all these algorithms, 
Gr\"obner bases are used as a fundamental tool.

\hspace{-3mm}{ \texttt{$\color{blue}>$ \color{red}ideal \color{black}AI =  \color{red}adjointIdeal\color{black}(f); }}

\hspace{-3mm}{ \texttt{$\color{blue}>$ AI; }}\color{blue}

\hspace{-3mm}{ \texttt{ \hspace{1mm} \_[1]=y3-y2z}}

\hspace{-3mm}{ \texttt{ \hspace{1mm} \_[2]=xy2-xyz}}

\hspace{-3mm}{ \texttt{ \hspace{1mm} \_[3]=x2y-xyz}}

\hspace{-3mm}{ \texttt{ \hspace{1mm} \_[4]=x3-x2z}}\color{black}

\noindent The resulting four cubic generators of the adjoint ideal define the curves depicted in 
Figure \ref{fig adjoints generators}, where the thickening of a line
indicates that the line comes with a double structure. A general adjoint curve, that 
is, a curve defined by a general linear combination of the four generators, 
is shown in Figure \ref{fig adjoint}.

\begin{figure}[h]

\begin{center}
\begin{tabular}
[c]{ll}%
\raisebox{0.0934in}{\includegraphics[
height=1.0309in,
width=1.6328in
]%
{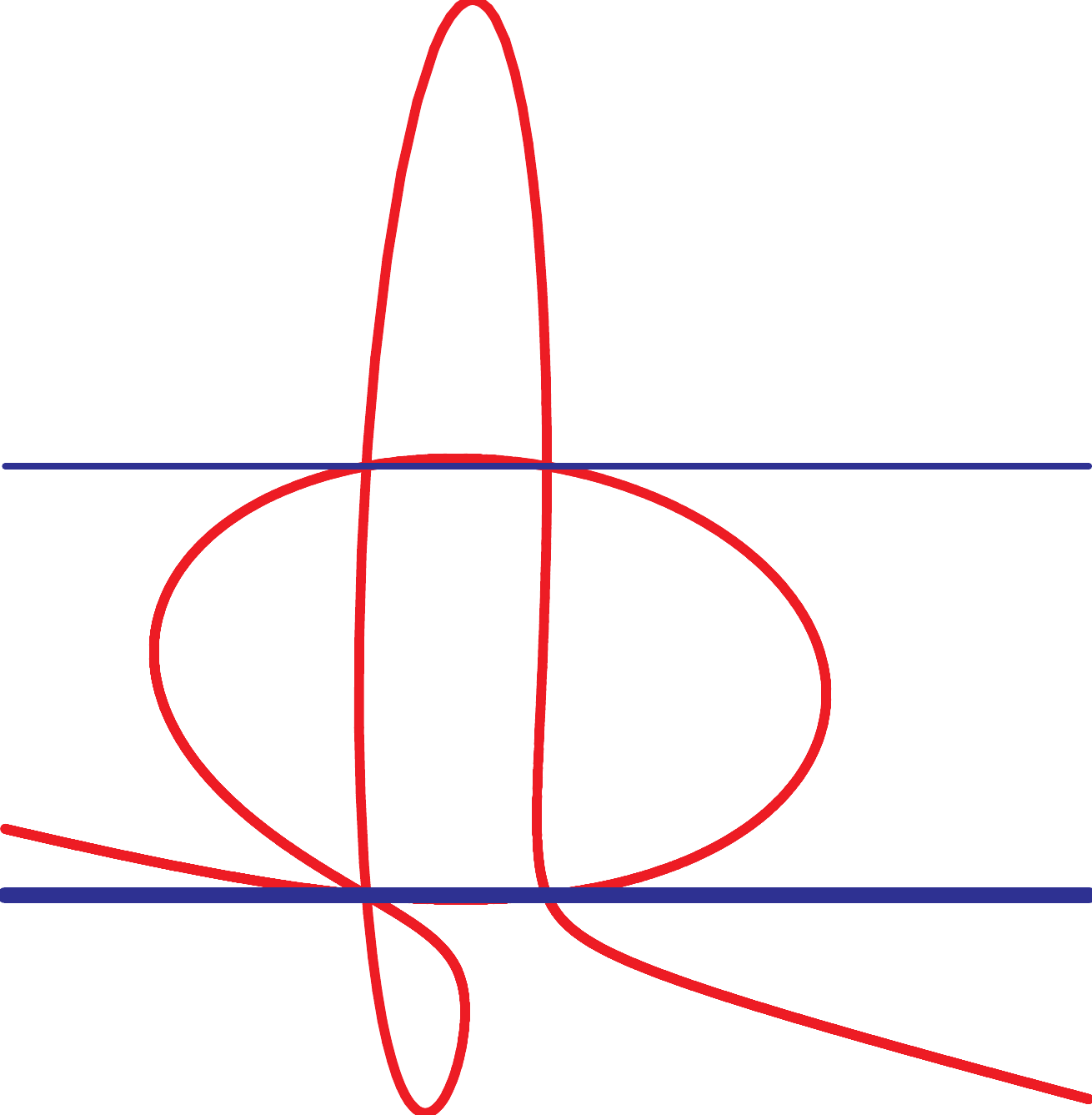}%
}
 &
{\includegraphics[
height=1.2324in,
width=1.6328in
]%
{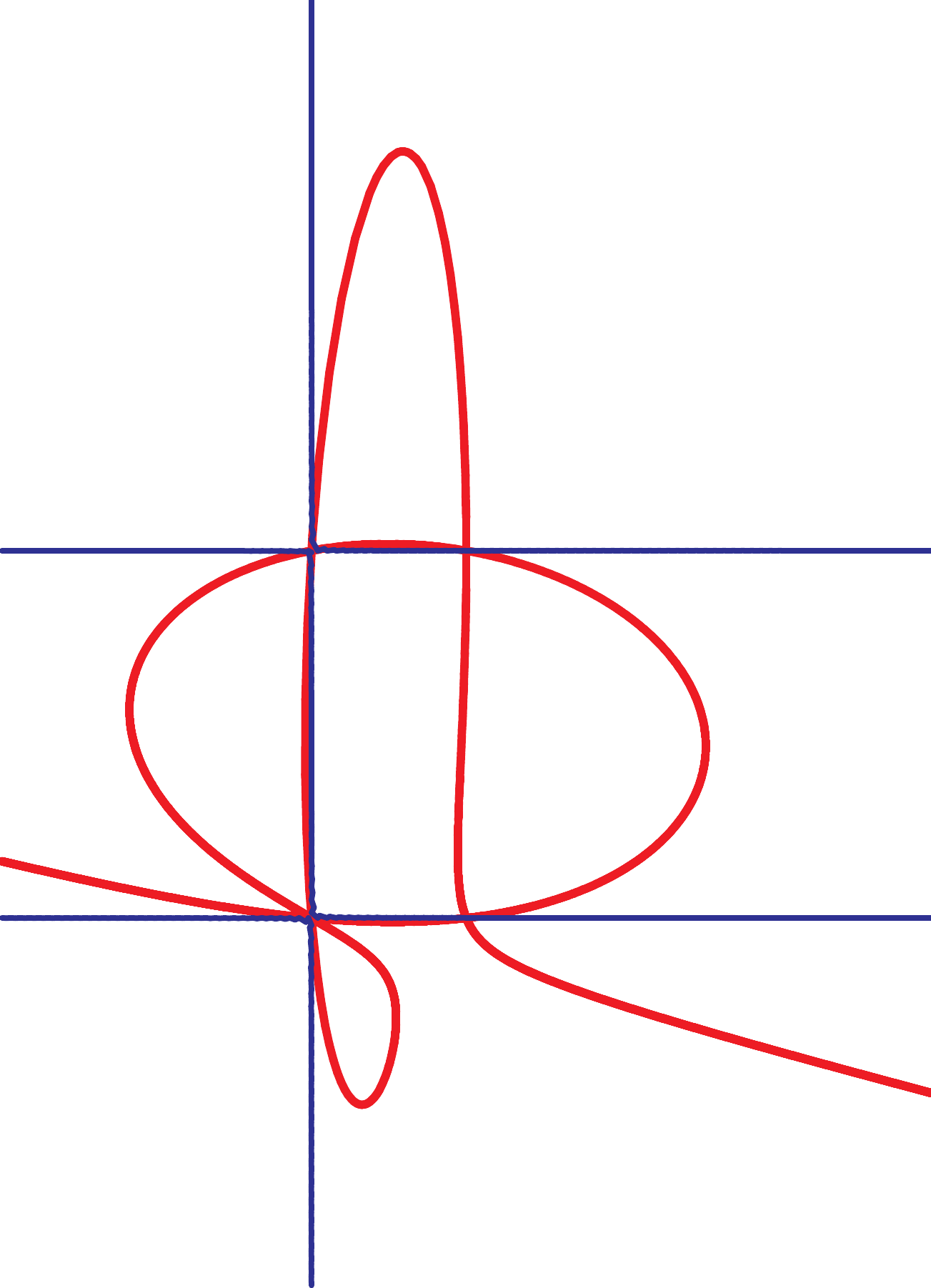}%
}
\\%
{\includegraphics[
height=1.2324in,
width=1.6336in
]%
{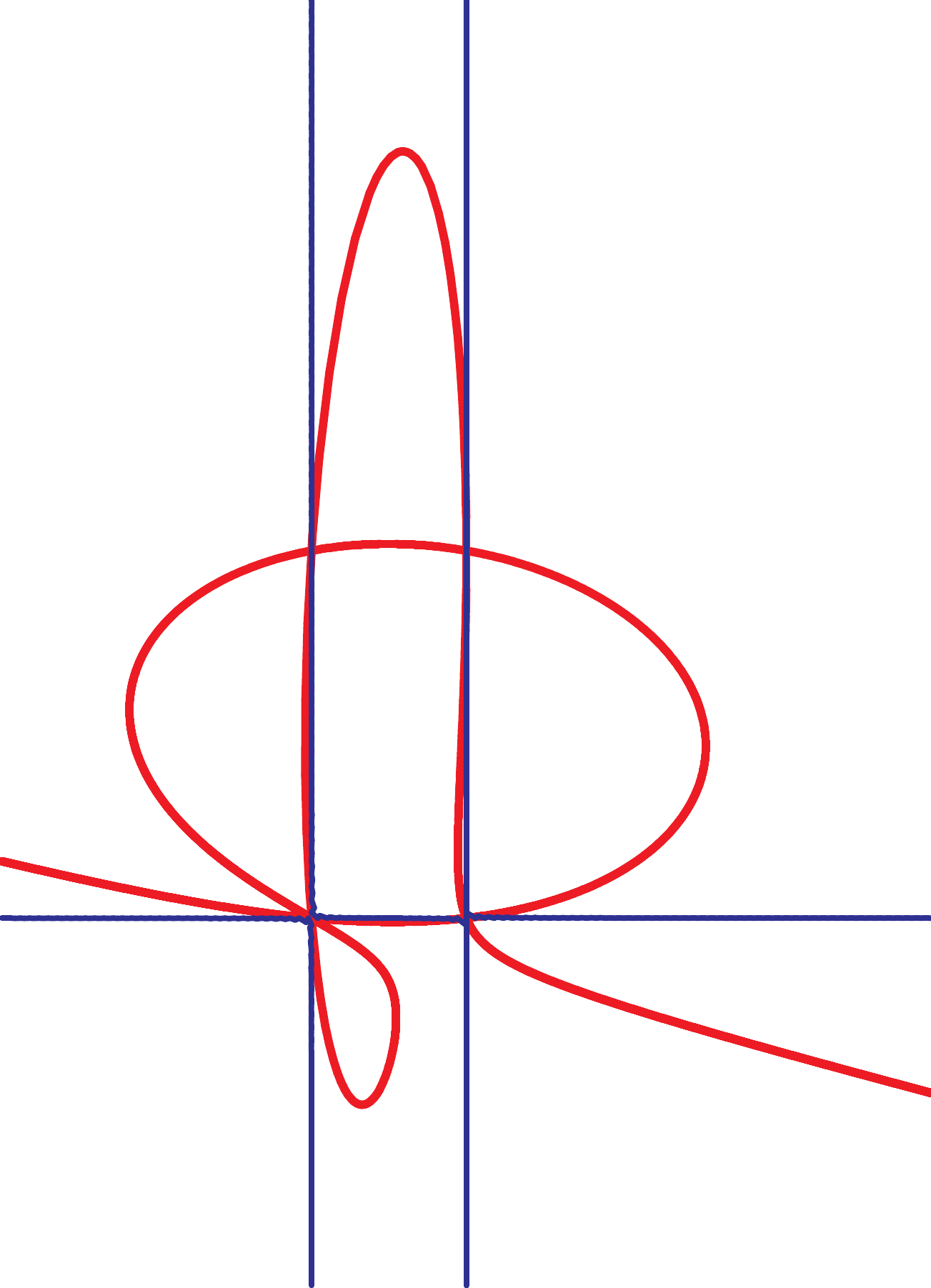}%
}
 &
{\includegraphics[
height=1.2324in,
width=1.6328in
]%
{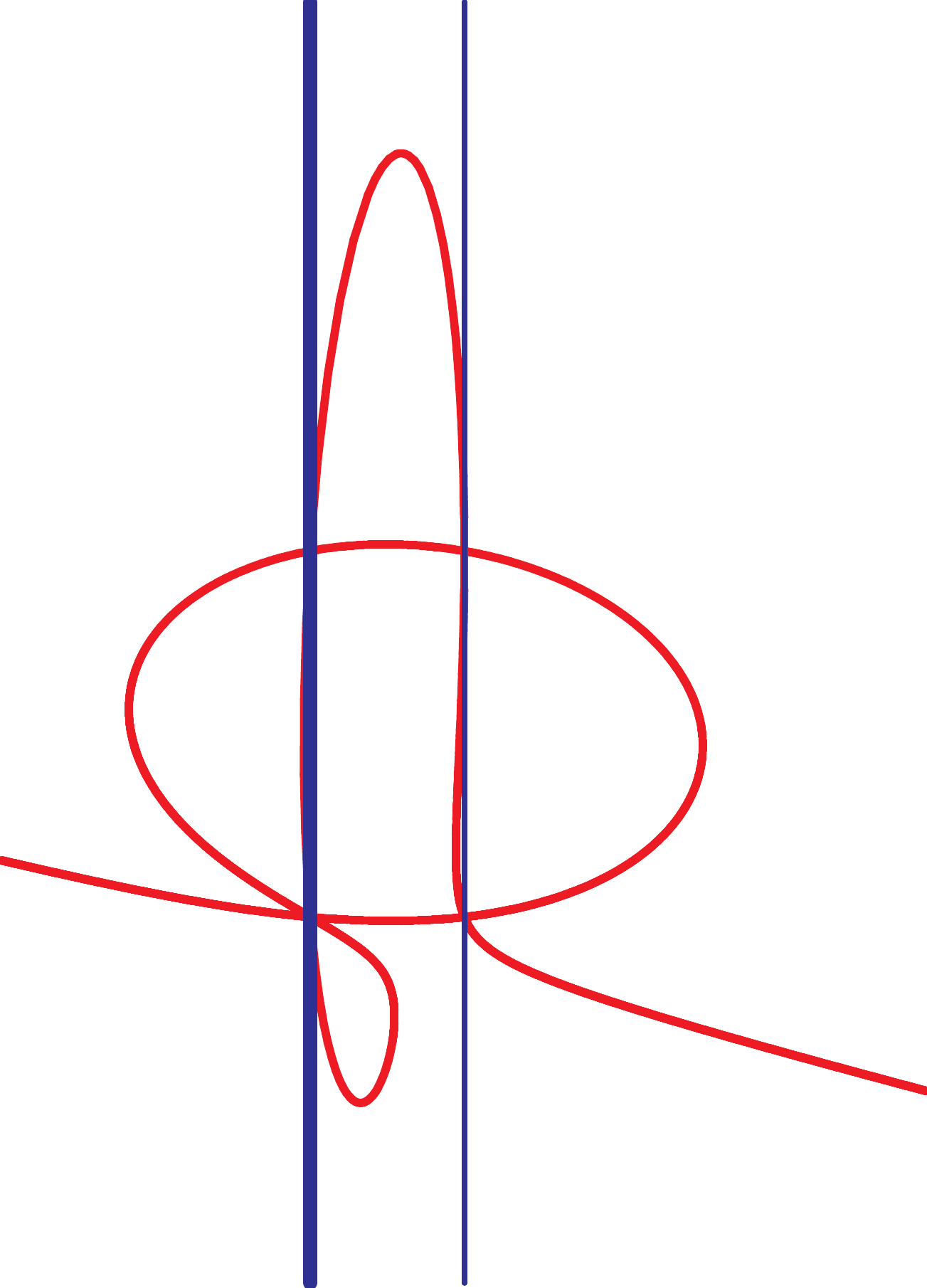}%
}
\end{tabular}
\end{center}
\caption{Cubic curves defined by the generators of the adjoint ideal of a degree-$5$ curve
with three ordinary double points and one ordinary triple point. The degree-5 curve
is shown in red.}%
\label{fig adjoints generators}%
\end{figure}

\begin{figure}[h]

\begin{center}
{\includegraphics[
height=2.46in,
width=3.265in
]%
{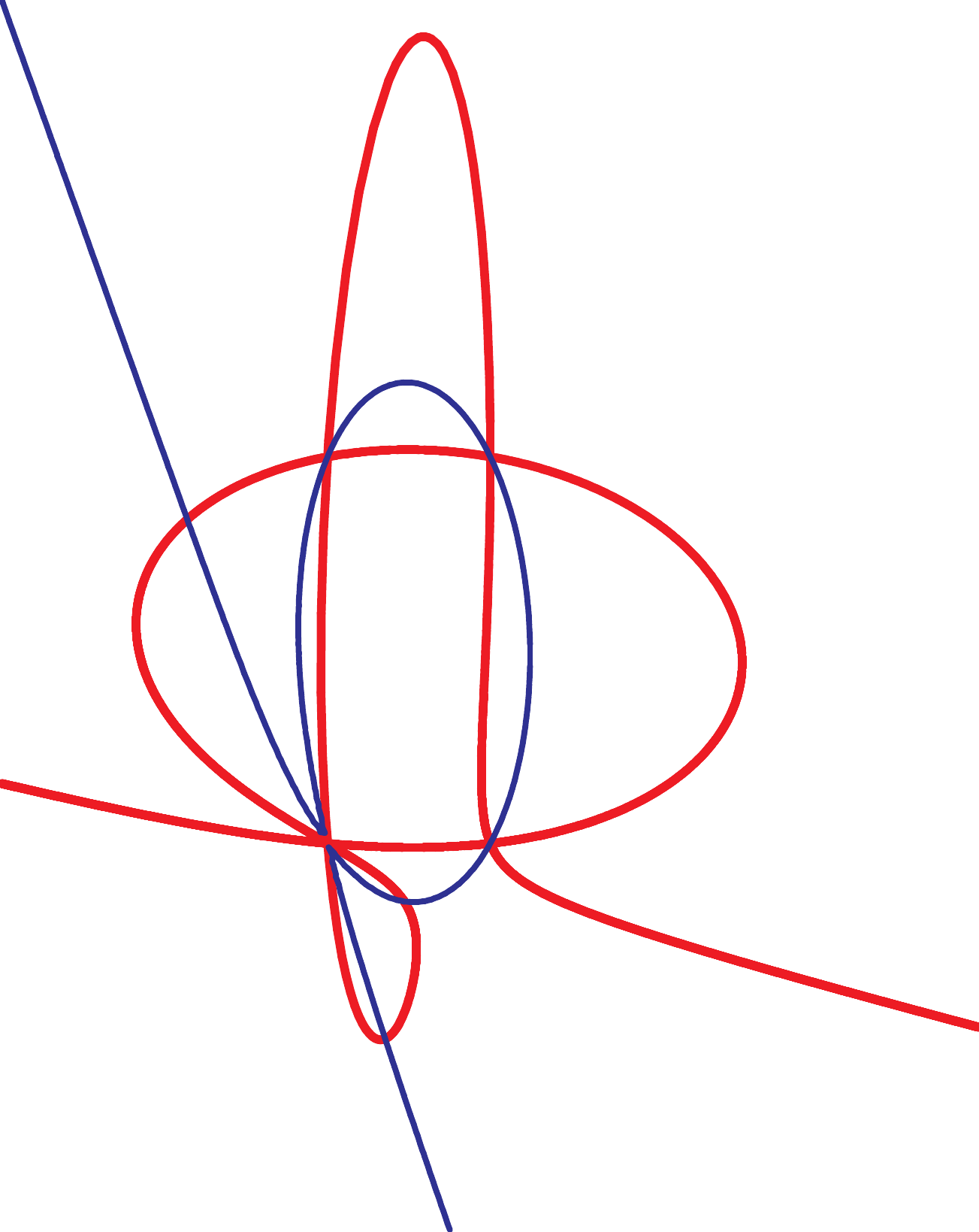}%
}
\end{center}
\caption{A general adjoint curve of $C$ of degree $3$.}%
\label{fig adjoint}%
\end{figure}

\noindent The four generators give a birational map from $C$  to a curve 
$\widetilde{C}$ in projective 3-space $\mathbb P^3$. We obtain $\widetilde{C}$ via elimination, 
a typical application of Gr\"obner bases:

\hspace{-3mm}{ \texttt{$\color{blue}>$ \color{red}def \color{black} Rn = \color{red}mapToRatNormCurve\color{black}(f,AI);}}

\hspace{-3mm}{ \texttt{$\color{blue}>$ \color{red}setring\color{black}(Rn);}}

\hspace{-3mm}{ \texttt{$\color{blue}>$ \color{black} RNC;}}\color{blue}

\hspace{-3mm}{ \texttt{ \hspace{1mm} RNC[1]=y(2)*y(3)-y(1)*y(4)}}

\hspace{-3mm}{ \texttt{ \hspace{1mm} RNC[2]=20*y(1)*y(2)-20*y(2)\symbol{94}2+130*y(1)*y(4)}}

\hspace{-3mm}{ \texttt{ \hspace{14mm}+20*y(2)*y(4)+10*y(3)*y(4)+y(4)\symbol{94}2}}

\hspace{-3mm}{ \texttt{ \hspace{1mm} RNC[3]=20*y(1)\symbol{94}2-20*y(1)*y(2)+130*y(1)*y(3)}}

\hspace{-3mm}{ \texttt{ \hspace{14mm}+10*y(3)\symbol{94}2+20*y(1)*y(4)+y(3)*y(4)}}\color{black}

\noindent Note that $\widetilde{C}$ is a variant of the projective twisted cubic curve, 
the rational normal curve in $\mathbb P^3$  (for a picture see Figure \ref{fig normalization}). 
This non-singular curve is mapped isomorphically onto the projective line 
$\mathbb{P}^1$ by the anticanonical linear system, which can be
computed using syzygies:

{ \texttt{$\color{blue}>$ \color{red}rncAntiCanonicalMap\color{black}(RNC);}}\color{blue}

\hspace{-3mm}{ \texttt{ \hspace{1mm} \_[1]=2*y(2)+13*y(4)}}

\hspace{-3mm}{ \texttt{ \hspace{1mm} \_[2]=y(4)}}\color{black}

\noindent 
Composing all maps in this construction, and inverting the resulting birational map,
we get the desired parametrization. 
In general, depending on the number of generators of the adjoint ideal, the rational
normal curve computed by the algorithm is embedded into a projective space of
odd or even dimension. In the latter case, successive applications of the canonical
linear system map the normal curve onto a plane conic. Computing a rational 
parametrization of the conic is equivalent to finding a point on the conic. It can 
be algorithmically decided whether we can find such a point with rational coordinates 
or not. In the latter case, we have to pass to a quadratic field extension of $\mathbb{Q}$.
\end{example}

\begin{remark} The need of passing to a field extension occurs in many geometric 
constructions. Often, repeated field extensions are needed. The effective 
computation of Gr\"obner bases over (towers of) number fields is therefore of 
utmost importance. One general way of achieving higher speed is the 
parallelization of algorithms. This will be addressed in the next section, where we 
will, in particular, discuss a parallel version of the Gr\"obner basis (syzygy) algorithm which
is specific to number fields \cite{derje}. New ideas for enhancing syzygy
computations in general are presented in \cite{SYZ}. Combining the two 
approaches in the case of number fields is a topic of future research.
\end{remark}

\subsection{Parallelization}

Parallelizing computer algebra systems allows for the efficient use of
multicore computers and high-performance clusters. To achieve parallelization is a tremendous 
challenge  both from a computer science and a mathematical point of view. 

From a computer science point of view, there are two possible approaches:

\begin{itemize}[leftmargin=*,itemsep=2pt]
\item  Distributed and multi-process systems work by using different processes that do not share memory and communicate 
by message passing. These systems only allow for \emph{coarse-grained parallelism}, which limits their ability to work on 
large shared data structures, but can in principle scale up indefinitely.

\item  Shared memory systems work by using multiple threads of control in a single process operating on shared data. 
They allow for more \emph{fine-grained parallelism} and more sophisticated concurrency control, down to the level of individual 
CPU instructions, but are limited in their scalability by how many processors can share efficient access to the same memory 
on current hardware.
\end{itemize}

For best performance, typically hybrid models are used, which exploit the strengths of both 
shared memory and distributed systems, while mitigating their respective downsides.

From its version 3.1.4 on,  \textsc{Singular} has been offering a framework for coarse-grained 
parallelization, with a convenient user access provided by the library \texttt{parallel.lib} 
\cite{parallellib}. The example below illustrates the use of this framework:

\begin{example}[Coarse  Grained Parallelization in \textsc{Singular}]
We implement a \textsc{Singular} procedure which computes a Gr\"obner basis for a given ideal with respect to a given
monomial ordering. The procedure returns  the size of the Gr\"obner basis. We apply it
in two parallel runs to a specific ideal in $\Q[x_1,\dots,x_4]$,
choosing for one run the lexicographical monomial ordering  
\texttt{lp} and for the other run the degree reverse lexicographical ordering~\texttt{dp}:

\hspace{-3mm}{ \texttt{$\color{blue}>$
\color{red}LIB \color{black}"parallel.lib"; \color{red}LIB \color{black}"random.lib";}}

\hspace{-3mm}{ \texttt{$\color{blue}>$ \color{red}proc
\color{black}sizeGb(\color{red}ideal \color{black}I, \color{red}string
\color{black}monord)\{} }

\hspace{-3mm}{ \qquad\qquad\texttt{\color{red}def \color{black}R =
\color{red}basering; \color{red}list \color{black}RL =
\color{red}ringlist\color{black}(R);}}

\hspace{-3mm}{ \qquad\qquad\texttt{RL[3][1][1] = monord; \color{red}def
\color{black}S = \color{red}ring\color{black}(RL);
\color{red}setring\color{black}(S);} }

\hspace{-3mm}{ \qquad\qquad
\texttt{\color{red}return\color{black}(\color{red}size\color{black}(\color{red}groebner\color{black}(\texttt{\color{red}imap\color{black}(R,I)}%
)));\}}}

\hspace{-3mm}{ \texttt{$\color{blue}>$ \color{red}ring \color{black}R =
0,x(1..4),dp;}}  

\hspace{-3mm}{ \texttt{$\color{blue}>$ \color{red}ideal \color{black}I =
\color{red}randomid\color{black}(\color{red}maxideal\color{black}(3),3,100);}}

\hspace{-3mm}{ \texttt{$\color{blue}>$ \color{red}list \color{black}commands =
"sizeGb","sizeGb";}}

\hspace{-3mm}{ \texttt{$\color{blue}>$ \color{red}list \color{black}args =
\color{red}list\color{black}(I,"lp"),\color{red}list\color{black}(I,"dp");}}

\hspace{-3mm}{ \texttt{$\color{blue}>$
\color{red}parallelWaitFirst\color{black}(commands, args);}}

\hspace{-3mm}{ \hspace{0.5cm}\texttt{\color{blue}[1] empty list}}

\hspace{-3mm}{ \hspace{0.5cm}\texttt{\color{blue}[2] 11}}

\hspace{-3mm}{ \texttt{$\color{blue}>$
\color{red}parallelWaitAll\color{black}(commands, args);}}

\hspace{-3mm}{ \hspace{0.5cm}\texttt{\color{blue}[1] 55}}

\hspace{-3mm}{ \hspace{0.5cm}\texttt{\color{blue}[2] 11}}

\noindent
As expected, the computation with respect to \texttt{dp} is much faster and leads to a Gr\"obner basis with less elements.
\end{example}
Using ideas from the successful parallelization of {\sc{GAP}} within the {\sc{HPC-GAP}} project
(see~\cite{reimer2,reimer3,reimer1}), a  multi-threaded prototype of {\sc{Singular}} has been 
implemented. Considerable further efforts are needed, however, to make this accessible to
users without a deep background in parallel programming.

From a mathematical point of view, there are algorithms whose basic strategy is inherently 
parallel, whereas others are sequential in nature.
A prominent example of the former type is Villamayor's constructive version of
Hironaka's desingularization theorem, which will be briefly discussed in Section \ref{sec hironaka}. 
A prominent example of the latter type is the classical Grauert-Remmert type algorithm for 
normalization, which will be addressed at some length in Section \ref{sec normal}.

The systematic design of parallel algorithms for applications which so far can only be handled
by sequential algorithms is a major task for the years to come. For normalization,
this problem has recently been solved~\cite{BDLSS}. Over the field of rational numbers, the new algorithm
becomes particularly powerful by combining it with modular methods,
see again Section \ref{sec normal}.

\begin{figure}[H]

\begin{center}
\begin{tikzpicture}[
  scale=0.75,
  transform shape,
  -latex,
  every node/.style={
    minimum height=0.8cm,
    minimum width=1cm,
  },
  box/.style={
    draw,
    anchor=north,
  },
  nobox/.style={
    anchor=north,
  },
]

\newlength{\leveldist}
\setlength{\leveldist}{1.62cm}

\begin{scope}[
  level distance=\leveldist,
  sibling distance=0.2cm,
  edge from parent/.style={
    draw,
    edge from parent path={
      let \p1=($(\tikzchildnode.north)+(0pt,2pt)$),
      \p2=($(\tikzparentnode.south)!0.95!(\p1)$)
      in (\tikzparentnode.south) -- (\x2,\y1)
    }
  }
]
\Tree
[ .\node [box] (I) {$\widetilde{I}$};
  [ .\node [box] (I1) {$\widetilde{I}_{p_1}$};
      \node [box] (Gp11) {$\widetilde{G}_{1,p_1}$};
      \edge [dashed];
      \node [nobox, minimum width=0cm, inner sep=0pt] (cdots21) {$\cdots$};
      \node [box] (Gp1r) {$\widetilde{G}_{r_{p_1},p_1}$};
  ]
  [ .\node [box] (I2) {$\widetilde{I}_{p_2}$};
      \node [box] (Gp21) {$\widetilde{G}_{1,p_2}$};
      \edge [dashed];
      \node [nobox, minimum width=0cm, inner sep=0pt] (cdots22) {$\cdots$};
      \node [box] (Gp2r) {$\widetilde{G}_{r_{p_2},p_2}$};
  ]
  \edge [dashed];
  [ .\node [nobox, minimum width=0cm] (cdots1) {$\cdots$};
      \edge [draw=none];
      \node [nobox, minimum width=0cm, inner sep=0pt] (codts23) {$\cdots$};
  ]
  [ .\node [box] (Ik) {$\widetilde{I}_{p_k}$};
      \node [box] (Gpk1) {$\widetilde{G}_{1,p_k}$};
      \edge [dashed];
      \node [nobox, minimum width=0cm, inner sep=0pt] (cdots24) {$\cdots$};
      \node [box] (Gpkr) {$\widetilde{G}_{r_{p_k},p_k}$};
  ]
]
\end{scope}

\node [box, below=2\leveldist of I1.north] (G1) {$\widetilde{G}_{p_1}$};
\node [box, below=2\leveldist of I2.north] (G2) {$\widetilde{G}_{p_2}$};
\node [nobox, below=2\leveldist of cdots1.north, minimum width=0pt] (cdots3)
  {$\cdots$};
\node [box, below=2\leveldist of Ik.north] (Gk) {$\widetilde{G}_{p_k}$};
\node [box, below=4\leveldist of I.north] (bottom)
  {Modular Reconstruction (over $\Q$)};

\newcommand{\arrow}[3]{\draw [#3] let \p1=($(#2.north)+(0pt,2pt)$),
      \p2=($(#1.south)!0.90!(\p1)$)
      in (#1.south) -- (\x2,\y1);}
\arrow{Gp11}{G1}{};
\arrow{cdots21}{G1}{dashed};
\arrow{Gp1r}{G1}{};
\arrow{Gp21}{G2}{};
\arrow{cdots22}{G2}{dashed};
\arrow{Gp2r}{G2}{};
\arrow{Gpk1}{Gk}{};
\arrow{cdots24}{Gk}{dashed};
\arrow{Gpkr}{Gk}{};
\draw (G1.south) -- ($(bottom.north)+(-9pt,2pt)$);
\draw (G2.south) -- ($(bottom.north)+(-2pt,2pt)$);
\draw [dashed] (cdots3.south) -- ($(bottom.north)+(2pt,2pt)$);
\draw (Gk.south) -- ($(bottom.north)+(11pt,2pt)$);

\node [right=0.2cm of Gpkr] (level2) {level 2};
\node at (level2|-I) {Input};
\node at (level2|-Ik) {level 1};
\node at (level2|-Gk) {level 3};
\node[anchor=west] at (cdots1.west|-I)
  {$\widetilde{I} \subset \Q[X,t],\; f \in \widetilde{I}$};
\end{tikzpicture}
\end{center}
\caption{Two-fold parallel modular approach to Gr\"obner bases over number fields.}%
\label{fig number fields}%
\end{figure}
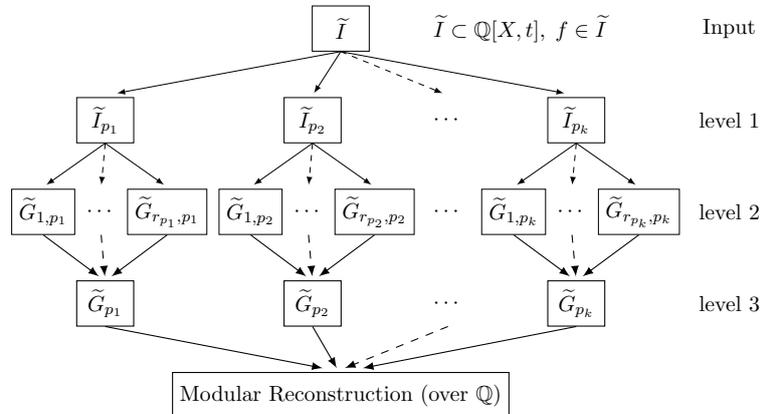

Modular methods are well-known for providing a way of parallelizing algorithms over 
$\mathbb Q$ (more generally, over number fields). For the fundamental task of
computing Gr\"obner bases, a modular version of Buchberger's 
algorithm is due to Arnold~\cite{arnold}. More recently, Boku, Fieker, 
Steenpa\ss\ and the second author~\cite{derje} have designed a modular
Gr\"obner bases algorithm which is specific to number fields. In addition 
to using the approach from Arnold's paper, which is to compute Gr\"obner 
bases modulo several primes and then use Chinese remaindering together 
with rational reconstruction, the new approach provides a second 
level of parallelization as depicted in Figure \ref{fig number fields}: If the
number field is presented as $K=\mathbb Q(\alpha)=\mathbb Q[t]/\langle f \rangle$, 
where $f\in \mathbb Q[t]$ is the minimal polynomial of $\alpha$, and if generators $g_1(X,\alpha),\dots, 
g_s(X,\alpha)$ for the ideal under consideration are given, represented by polynomials
$g_1(X,t),\dots, g_s(X,t)\in \mathbb Q[X,t]=Q[x_1,\dots,x_n,t]$, we wish to compute
a Gr\"obner basis for the ideal $\widetilde{I}=\langle g_1(X,t),\dots, g_s(X,t), f\rangle 
\subset \mathbb Q[X,t]$. The idea then is to reduce $\widetilde{I}$ modulo a suitable
number of primes $p_1,\dots, p_k$ (level 1 of the algorithm), get the second level of parallelization
by factorizing the reductions of $f$ modulo the $p_i$, and use, for each $i$,
polynomial Chinese remaindering to put the results modulo $p_i$ together
(level 3 of the algorithm).

\subsection{Make More and More of the Abstract Concepts 
of Algebraic Geometry Constructive}\label{sec hironaka}

The following groundbreaking theorem proved by Hironaka in $1964$ shows the existence of resolutions of
singularites in characteristic zero. It is worth mentioning that, on his way, Hironaka introduced the idea of 
standard bases, the power series analogue of Gr\"obner bases.

\begin{theorem}
[Hironaka, 1964]For every algebraic variety over a field $K$ of characteristic zero,
a desingularization can be obtained by a finite sequence of blow-ups along smooth centers.
\end{theorem}

We illustrate the blow-up process by a simple example:

\begin{example}
As shown in Figure \ref{node}, a node can be resolved by a single blow-up: we replace the node 
by a line and separate, thus,  the two branches of the curve intersecting in the singularity.

\begin{figure}[h]

\begin{center}
\includegraphics[
height=1.6646in,
width=1.6646in
]%
{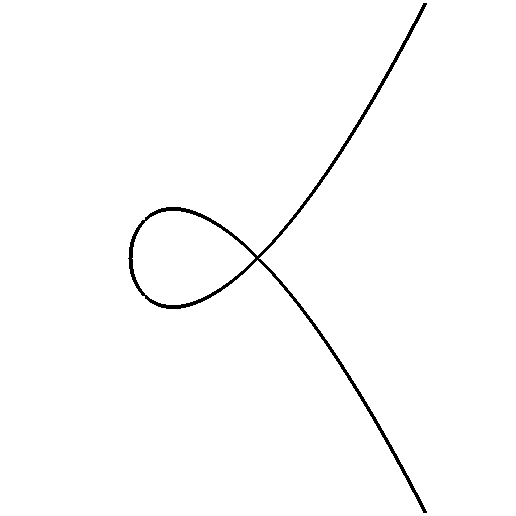}
\hspace{1cm}
\includegraphics[
height=1.8506in,
width=1.369in
]%
{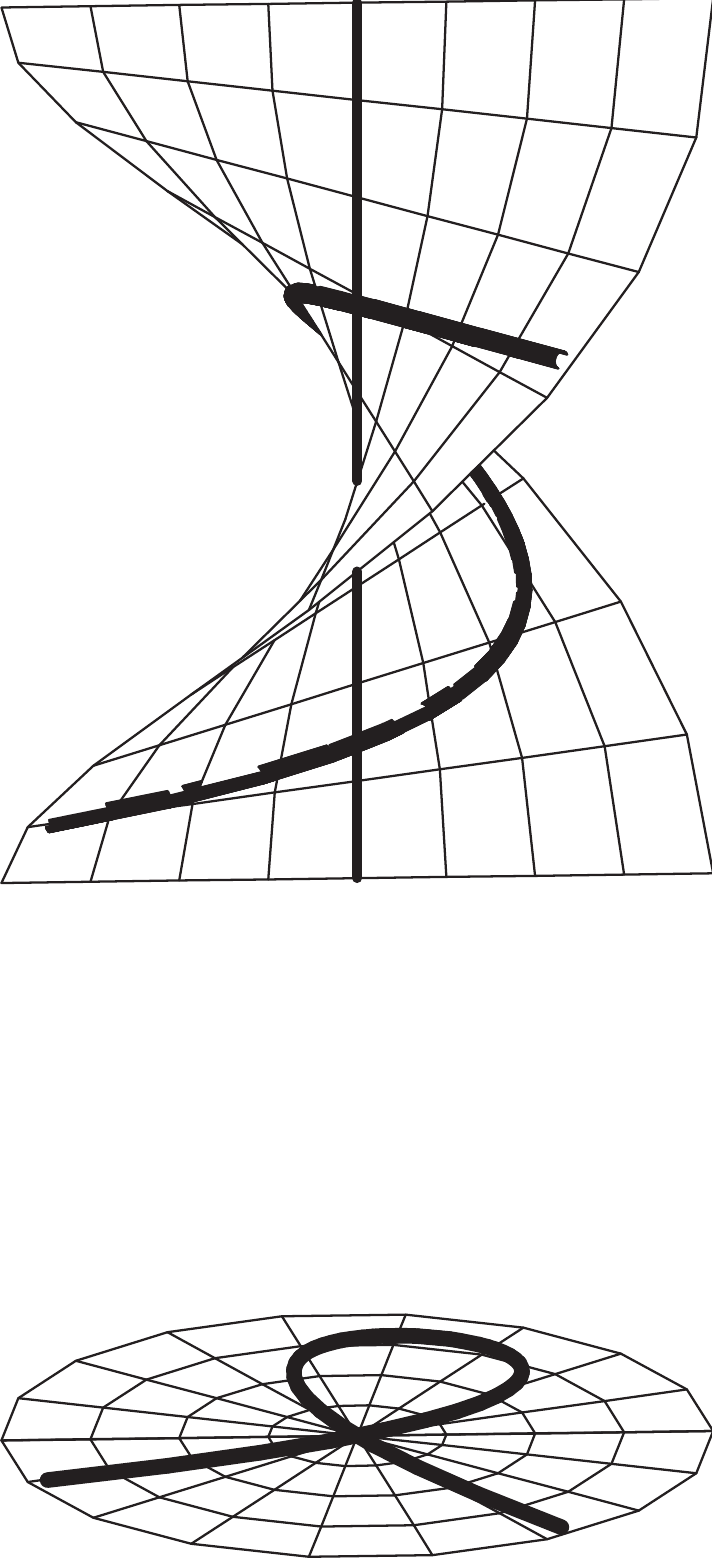}%

\end{center}
\caption{Blowing up a node.}%
\label{node}%
\end{figure}

\end{example}

In~\cite{BM,BEV,EH,FK1}, the abstract concepts developed by Hironaka have been translated 
into an algorithmic approach to desingularization. An effective variant of this, 
relying on a clever selection of the centers for the blow-ups, has been 
implemented by Fr{\"u}hbis-Kr{\"u}ger and Pfister in the \textsc{Singular} library 
\texttt{resolve.lib}~\cite{resolve}.

\begin{figure}[h]

\begin{center}
\includegraphics[
width=2in
]%
{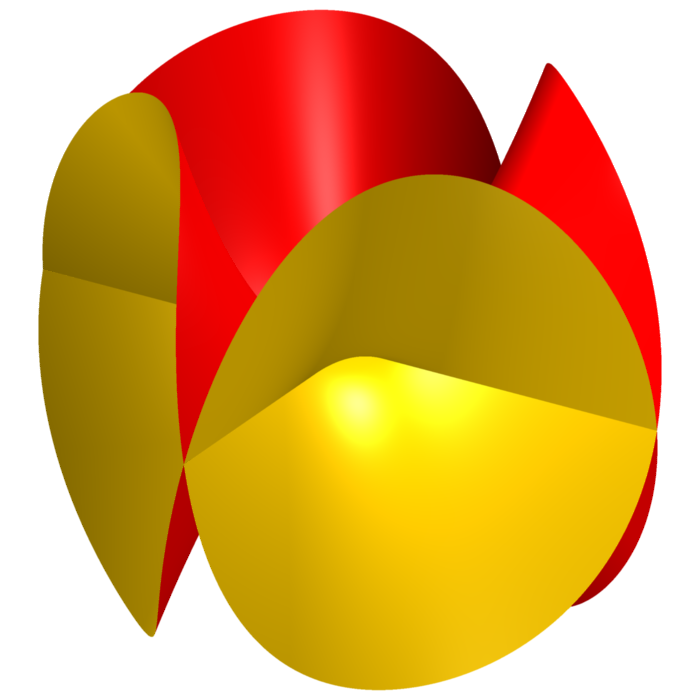}%
\end{center}
\caption{The surface $z^{2}-x^{2} y^{2}=0$.}%
\label{surface}%
\end{figure}

\begin{figure}[h]

\begin{center}
\includegraphics[
width=4in
]%
{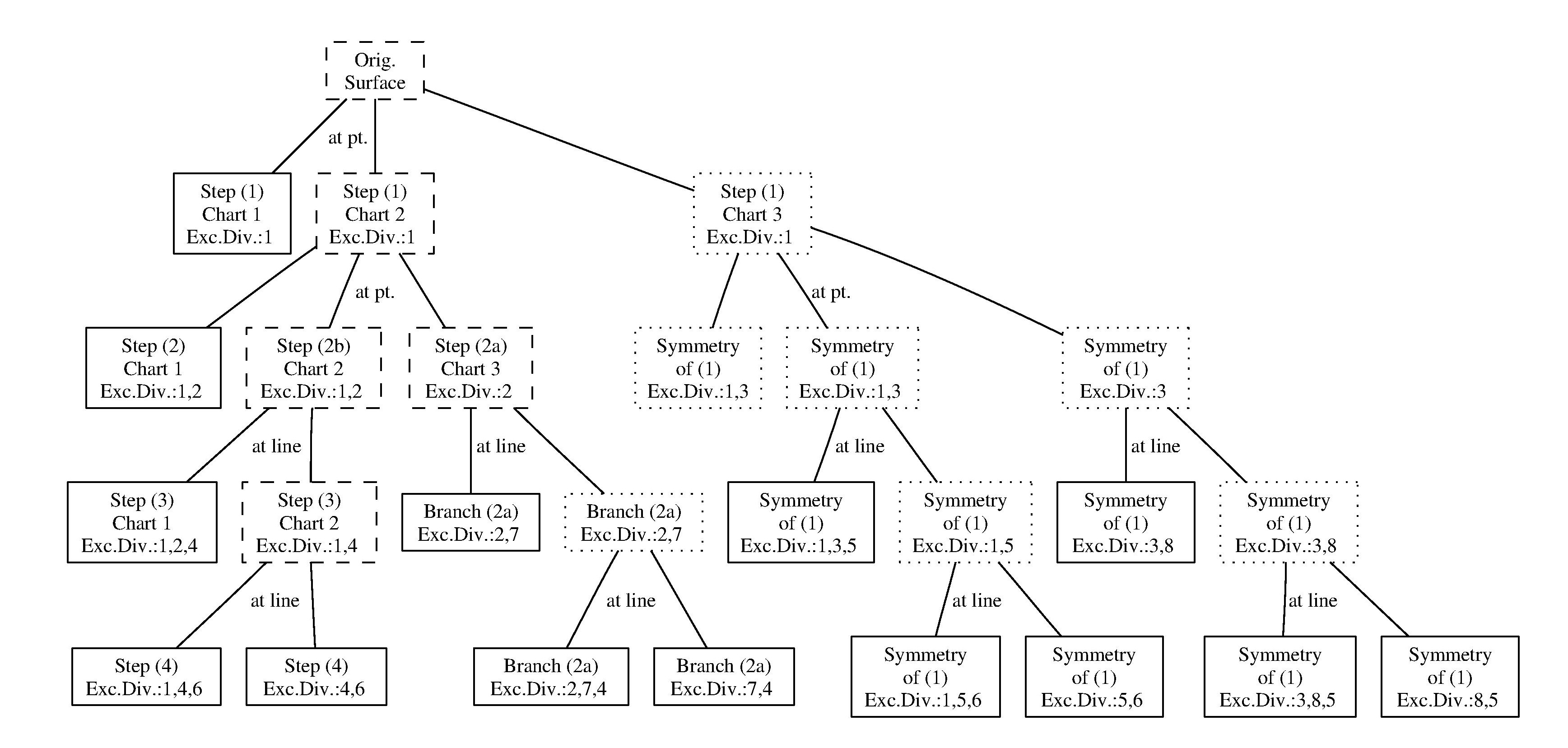}%
\end{center}
\caption{The tree of charts.}%
\label{resolve}%
\end{figure}

The desingularization algorithm is parallel in nature: Working with blow-ups means to work with 
different charts of projective spaces. In this way, the resolution of singularities 
leads to a tree of charts. Figure \ref{resolve} shows the graph for resolving the 
singularities of the hypersurface $z^{2}-x^{2}y^{2}=0$ which, in turn, is depicted in 
Figure \ref{surface}.%

Making abstract concepts constructive allows for both a better understanding 
of deep mathematical results and a computational treatment of the concepts.
A further preeminent example for this is the constructive version of the  Bernstein-Gel'fand-Gel'fand 
correspondence (BGG-correspondence)
by Eisenbud, Fl\o ystad, and Schreyer~\cite{zbMATH01963988}. This allows one
to express properties of sheaves over projective spaces in terms of exterior algebras.
More precisely, if $\mathbb P(V)$ is the projective space of lines in a vector space 
$V$, and $E$ is the exterior algebra $E=\Lambda V$, then the BGG-correspondence
relates coherent sheaves over $\mathbb P(V)$ to free resolutions over $E$. 
Since $E$ contains only finitely many monomials,  (non-commutative) Gr\"obner basis 
and syzygy computations over $E$ are often preferable to (commutative) Gr\"obner basis 
and syzygy computations over the homogeneous coordinate ring of $\mathbb P(V)$. 
One striking application of this, which is 
implemented in \textsc{Macaulay2}~\cite{GS} and \textsc{Singular}, gives a fast way 
of computing sheaf cohomology. Providing computational access to cohomology
in all its disguises is a long-term goal of computational algebraic geometry.

The BGG-correspondence is an example of an equivalence of derived categories.
As we can see from the above discussion, such equivalences are not only interesting 
from a theoretical point of view, but may also allow for creating more effective algorithms --
provided they can be accessed computationally.

\subsection{Interaction and Integration of Computer Algebra Systems
and Libraries From Different Areas of Research}\label{sec systems}

\begin{figure}[h]

\begin{center}
{\includegraphics[
width=2.1353in
]%
{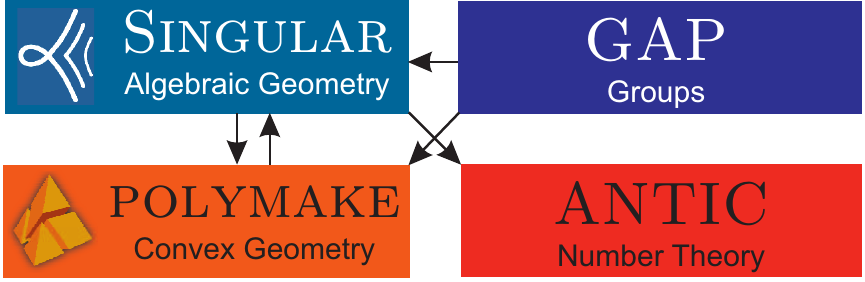}%
}

\end{center}
\caption{Directed interfaces.}%
\label{fig systems}%
\end{figure}

On the theoretical side, mathematical breakthroughs are often obtained by combining methods from 
different areas of mathematics. Making such connections accessible to computational methods
is another major challenge. Handling this challenge requires, in particular, that computer algebra
systems specializing in different areas are connected in a suitable way. One goal of the
Priority Programme SPP 1489, which was already mentioned in the introduction,  
is to interconnect \textsc{GAP}, \textsc{polymake}, \textsc{Singular}, and \textsc{Antic}.
So far, this has lead to directed interfaces as indicated in Figure \ref{fig systems}, with further 
directions and a much tighter integration of the systems subject to future development.

\begin{figure}[h]
\begin{center}
    \begin{tikzpicture}
      \node[anchor=east] (algCubic) at (0,0) {\includegraphics[scale=0.04]{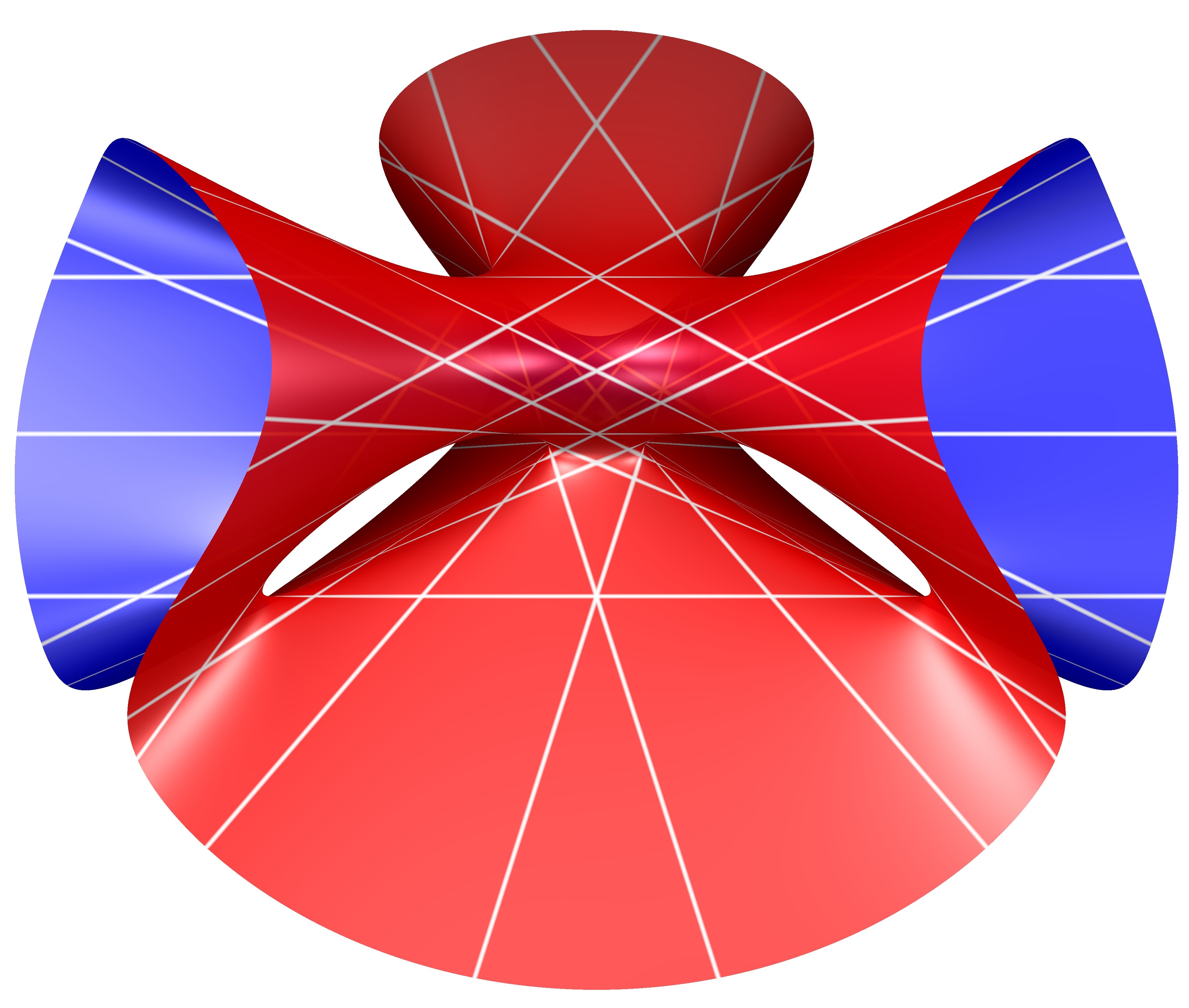}};
      \draw[very thick,->] (0,0) -- node[above,font=] {tropicalization} (3,0);
      \node[anchor=west] (tropCubic) at (3,0) {\includegraphics[scale=0.12]{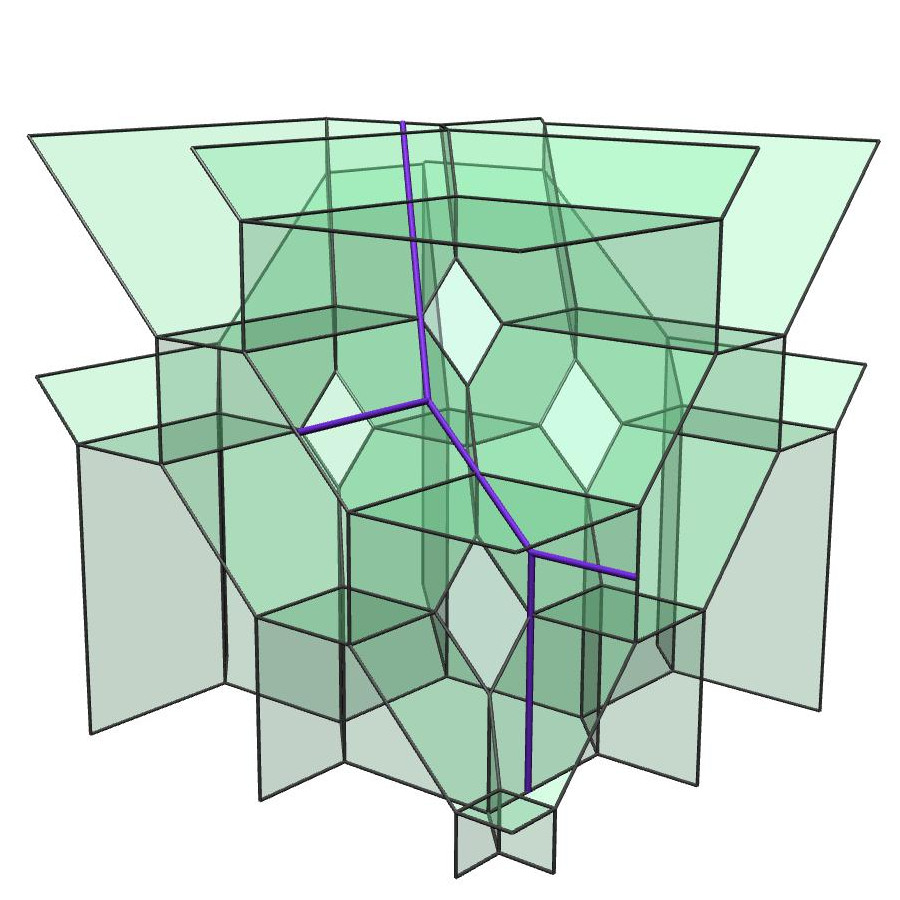}};
    \end{tikzpicture}
  \end{center}
\caption{The Tropicalization of Clebsch's diagonal cubic.}%
\label{fig:tropicalization}%
\end{figure}

In fact, the picture is much more complicated: The four systems rely on further
systems and libraries such as \textsc{normaliz}~\cite{normaliz} (affine monoids) 
and \textsc{Flint}~\cite{flint} (number theory),
and there are other packages which use at least one of the four systems, for example
\texttt{homalg}~\cite{homalg} (homological algebra) and \texttt{a-tint}~\cite{atint} 
(tropical intersection theory) (Figure \ref{fig:tropicalization}).

With regard to mathematical applications, the value of connecting \textsc{GAP} and \textsc{Singular}
is nicely demonstrated by Barakat's work on a several years old question of Serre to find a prediction 
for the number of connected components of unitary groups of group algebras in characteristic 2~\cite{serre,barakat}. 

A showcase application for combining \textsc{Singular}, \textsc{polymake}, and \textsc{GAP} 
is the symmetric algorithm for computing GIT-fans~\cite{BKR} by the first, third and fourth author~\cite{BKR}.
This algorithm, which will be discussed in more detail in Section \ref{sec git}, combines Gr\"obner basis and 
convex hull computations, and can make use of actions of finite symmetry groups.

\subsection{A Convenient Hierarchy of Languages}

Most modern computer algebra systems consist of two major components, 
a kernel which is typically written in C/C++ and a high level language for direct
user interaction, which in particular provides a convenient way for users to 
extend the system. While the kernel code is precompiled and, thus, performant,
the user language is interpreted, which means that it operates at a significantly 
slower speed. In addition to the differences in speed, the languages involved
provide different levels of abstraction with regard to modeling mathematical 
concepts. In view of the integration of different systems, a number of languages has
to be considered, leading to an even more complicated situation. To achieve 
the required level of performance and abstraction in this context, we need to set
up a convenient hierarchy of languages. Here, we propose in particular to examine 
the use of just-in-time compiled languages such as {\sc{Julia}}.

\subsection{Create and Integrate Electronic Libraries and Databases Relevant to Research}

Electronic libraries and databases of certain classes of mathematical
objects provide extremely useful tools for research in their respective fields. 
An example from group theory is the \emph{SmallGroups} library, which is 
distributed as a  {\sc{GAP}} package. An example from algebraic geometry is 
the \emph{Graded Ring Database},\footnote{See \url{http://www.grdb.co.uk}}
written by Gavin Brown and Alexander Kasprzyk, with contributions by several
other authors. The creation of such databases often depends on several computer
algebra systems. On the other hand, a researcher using the data may wish to 
access the database within a system with which he is already familiar.
This illustrates the benefits of a standardized approach to connect computer
algebra systems and mathematical databases.

\subsection{Facilitating the Access to Computer Algebra Systems}

\begin{figure}[h]
\begin{center}
\includegraphics[width = .9\textwidth]{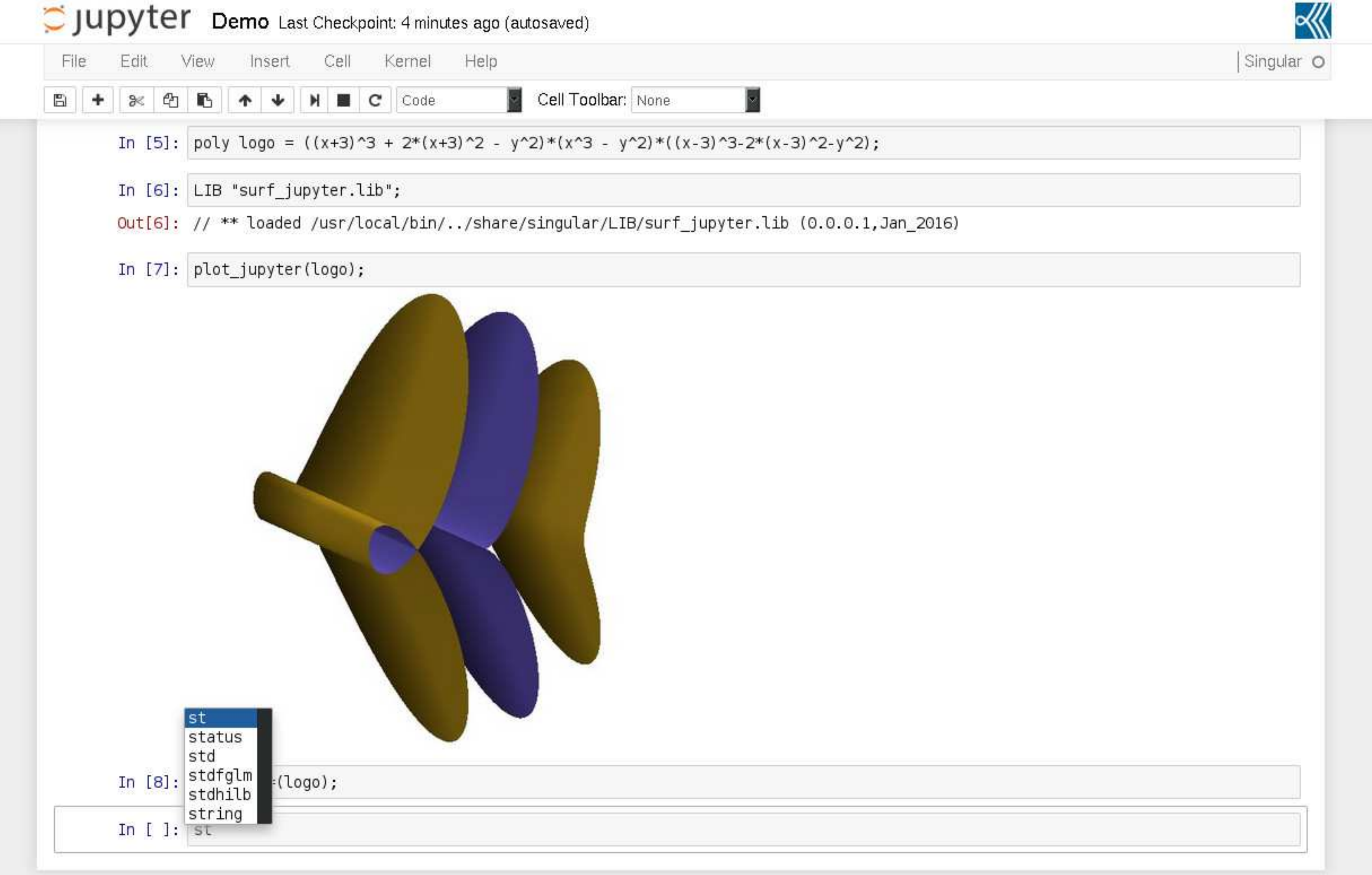}
\end{center}
\caption{Jupyter notebook for \textsc{Singular}.}%
\label{fig jupyter}%
\end{figure}

Computational algebraic geometry (and computer algebra in general) has a rapidly increasing 
amount of applications outside its original core areas, for example to computational biology, 
algebraic vision, and physics. As more and more non-specialists wish to use computer
algebra systems, the question of how to considerably ease the access to the systems arises
also in the Open Source community.
Virtual research environments such as the one developed within the \emph{OpenDreamKit}
project\footnote{See \url{http://opendreamkit.org}} may provide an answer to this 
question. Creating Jupyter notebooks\footnote{See \url{http://jupyter.org}} 
for systems such as {\sc{GAP}} and {\sc{Singular}} is one of the many goals of this project.
A  {\sc{Singular}} prototype has been written by Sebastian Gutsche, see Figure~\ref{fig jupyter}.

\section{A Parallel Approach to Normalization}\label{sec normal}

In this section, focusing on the normalization of rings, we give an example 
of how ideas from commutative algebra can be used to turn a
sequential algorithm into a parallel algorithm.

The normalization of rings is an important concept in commutative algebra, with
applications in algebraic geometry and singularity theory. Geometrically,
normalization removes singularities in codimension one and ``improves''
singularities in higher codimension. In particular, for curves, normalization yields a 
desingularization (see Examples \ref{ex normalization} and \ref{ex:normal} below). 
From a computer algebra point of view, normalization is fundamental to quite a 
number of algorithms with applications in algebra, geometry, and number theory.
In Example \ref{ex:para}, for instance, we have used normalization 
to compute adjoint curves and, thus, parametrizations of  rational curves.

The by now classical Grauert-Remmert type approach~\cite{deJong98,DGPJ,GLS} 
to compute normalization proceeds by successively enlarging the given 
ring until the Grauert-Remmert normality criterion~\cite{GR}  tells us that the 
normalization has been reached. Obviously, this approach is completely sequential 
in nature. As already pointed out, it is a major challenge to systematically 
design parallel alternatives to basic and high-level algorithms which are 
sequential in nature. For normalization, this problem has recently been  solved 
in~\cite{BDLSS} by using the technique of localization and proving a local version
of the Grauert-Remmert normality criterion.

To explain this in more detail, we suppose for simplicity
that the ring under consideration is an affine domain over a field $K$. That is, 
we consider a quotient ring of type 
$A=K[x_{1},\ldots,x_{n}]/I$, where $I$ is a prime ideal.
We require that $K$ is a perfect field.

We begin by recalling some basic definitions and results.

\begin{definition}
The \textbf{normalization}  of $A$ is the integral closure $\overline{A}$
of $A$ in its quotient field $Q(A)$,
$$
\overline{A} =\{a\in Q(A) \mid \text{there exists } f\in A[t]\text{ monic with }f(a)=0\}.
$$
We call $A$ \textbf{normal} if
$A=\overline{A}$.
\end{definition}

By Emmy Noether's finiteness theorem (see~\cite{HS}), we may represent 
$\overline{A}$ as the set of $A$-linear combinations of a finite set of elements of
$\overline{A}$. That is:

\begin{theorem}
[Emmy Noether]\label{thm:noether}$\overline{A}$ is a finitely generated $A$-module.
\end{theorem}

We also say that the ring extension $A\subset \overline{A}$ is finite.
In particular, $\overline{A}$  is again an affine domain over $K$.

\begin{example}
\label{ex normalization}For the coordinate ring $A=K[x,y]/I$ of the 
nodal plane curve $C=V(I)$ defined by the prime ideal 
$I=\left\langle x^{3}+x^{2}-y^{2}\right\rangle \subset K[x,y]$,
we have%
\[%
\begin{tabular}
[c]{rcccc}%
$A=K[x,y]/I$ & $\cong$ & $K[t^{2}-1,t^{3}-t]$ & $\subset$ & $K[t]\cong%
\overline{A}$.\\
$\overline{x}$ & $\mapsto$ & \multicolumn{1}{l}{$t^{2}-1$} &  & \\
$\overline{y}$ & $\mapsto$ & \multicolumn{1}{l}{$t^{3}-t$}&  &
\end{tabular}
\]
In particular, $\overline{A}$ is generated as an $A$-module by $1$ and $\frac
{\overline{y}}{\overline{x}}$. 

Geometrically, the inclusion map $A\hookrightarrow\overline{A}$
corresponds to the parametrization%
\begin{align*}
\mathbb A^1(K)  & \rightarrow C\subset \mathbb A^2(K),\qquad
t \mapsto(t^{2}-1,t^{3}-t).
\end{align*}
In other words,  the parametrization is the normalization (desingularization) 
map of  the rational curve $C$.
\end{example}

Historically, the first Grauert-Remmert-type algorithm for normalization is due to 
de Jong~\cite{deJong98,DGPJ}. This algorithm has been implemented in \textsc{Singular}, 
\textsc{Macaulay2}, and \textsc{{Magma}}~\cite{BCP}.
The algorithm of Greuel, Laplagne, and Seelisch~\cite{GLS} is a more efficient
version of de Jong's algorithm. It is implemented in the \textsc{{Singular}}
library \texttt{normal.lib}~\cite{normal}. 

The starting point of these algorithms is the following lemma:

\begin{lemma}[\cite{GLS}]
If $J\subset A$ is an ideal and $0\neq g\in J$, then there are natural inclusions of rings
\[
A\, \hookrightarrow\, \operatorname*{Hom}\nolimits_{A} (J, J)\,\cong\,\frac{1}{g}(gJ :_{A} J)\,\subset\,\overline{A}%
\,\subset\,{\text{Q}}(A),\quad a \mapsto\varphi_{a}, \quad \varphi\mapsto\frac
{\varphi(g)}{g},
\]
where $\varphi_{a}$ is the multiplication by $a$.
\end{lemma}

Now, starting from $A_{0}=A$
and $J_{0}=J$, and setting%
\[%
\begin{tabular}
[c]{lll}%
$A_{i+1}=\frac{1}{g}(gJ_{i}:_{A_{i}}J_{i})$ & and  & $J_{i}=\sqrt{JA_{i}}$,%
\end{tabular}
\
\]
we get a chain of finite extensions of affine domains which becomes eventually stationary 
by Theorem \ref{thm:noether}:
\[
A=A_{0}\subset\dots\subset A_{i}\subset\dots\subset A_{m}=A_{m+1} \subset \overline{A}.
\]
The Grauert-Remmert-criterion for normality tells us that for an appropriate
choice of $J$, the process described above terminates with the normalization $A_{m}=\overline{A}$.
In formulating the criterion, we write $N(A)$ for the \textbf{non-normal locus} 
of $A$, that is, if
$$
\operatorname{Spec}(A)=\{P\subset A\mid P{\text{ prime ideal}}\}
$$
denotes the \textbf{spectrum} of $A$, and $A_P$ the localization of $A$ at $P$, then
\[
N(A)=\{P\in\operatorname{Spec}(A)\mid A_{P}\text{ is not normal}\}.
\]
\begin{theorem}
[Grauert-Remmert~\cite{GR}] Let $\langle 0\rangle\neq J\subset A$ be an ideal with $J=\sqrt
{J}$ and such that
\[
N(A)\subseteq V(J):=\{P \in \operatorname{Spec}(A)\mid P\supset J\}\text{.}%
\]
Then $A$ is normal if and only if $A \cong\operatorname*{Hom}\nolimits_{A}(J, J)$ via the map
which sends $a$ to multiplication by $a$.
\end{theorem}
The problem now is that we do not know an algorithm for computing $N(A)$,
except if the normalization is already known to us. To remedy this situation,
we consider the \textbf{singular locus} of $A$,
\[
\operatorname{Sing}(A)=\{P\in\operatorname{Spec}(A)\mid A_{P}\text{ is not
regular}\},
\]
which contains the non-normal locus:  $N(A)\subseteq \operatorname{Sing}(A)$. 
Since we work over a perfect field $K$, the Jacobian
criterion tells us that $\operatorname{Sing}(A) = V(\operatorname{Jac}(I))$,
where $\operatorname{Jac}(I)$ is the Jacobian ideal\footnote{The 
Jacobian ideal of $A$ is generated by the images of the $c\times c$ minors of the 
Jacobian matrix $(\frac{\partial f_i}{\partial x_j})$, where $c$ is the codimension 
and $f_1,\dots,  f_r$ are polynomial generators for $I$.} 
of $A$ (see~\cite{Eis}).
Hence, if we choose $J=\sqrt{\operatorname{Jac}(I)}$, the above process terminates
with $A_{m}=\overline{A}$ by the following lemma.

\begin{lemma}[\cite{GLS}]
With notation as above, $N(A_{i})\subseteq V(\sqrt{JA_{i}})$ for all $i$.
\end{lemma}

\begin{example}\label{ex:normal}
For the coordinate ring $A$ of the plane algebraic curve $C$ from Example
\ref{ex normalization}, the normalization algorithm returns the coordinate
ring of a variant of the twisted cubic curve $\overline C$ in affine $3$-space, where the
inclusion $A\subset\overline{A}$ corresponds to the projection of
$\overline C$ to $C$ via $(x,y,z)\mapsto (x,y)$ as shown in Figure \ref{fig normalization}.
This result fits with the result in Example \ref{ex normalization}:
The curve  $\overline C$ is rational, with a parametrization given by
\begin{align*}
\mathbb A^1(K)  & \rightarrow \overline{C}\subset \mathbb A^3(K),\qquad
t \mapsto(t^{2}-1,t^{3}-t,t).
\end{align*}
Composing this with the projection, we get the normalization map from 
Example~\ref{ex normalization}.
\begin{figure}[H]
\begin{center}
\includegraphics[
height=2.559in,
width=2.0998in
]%
{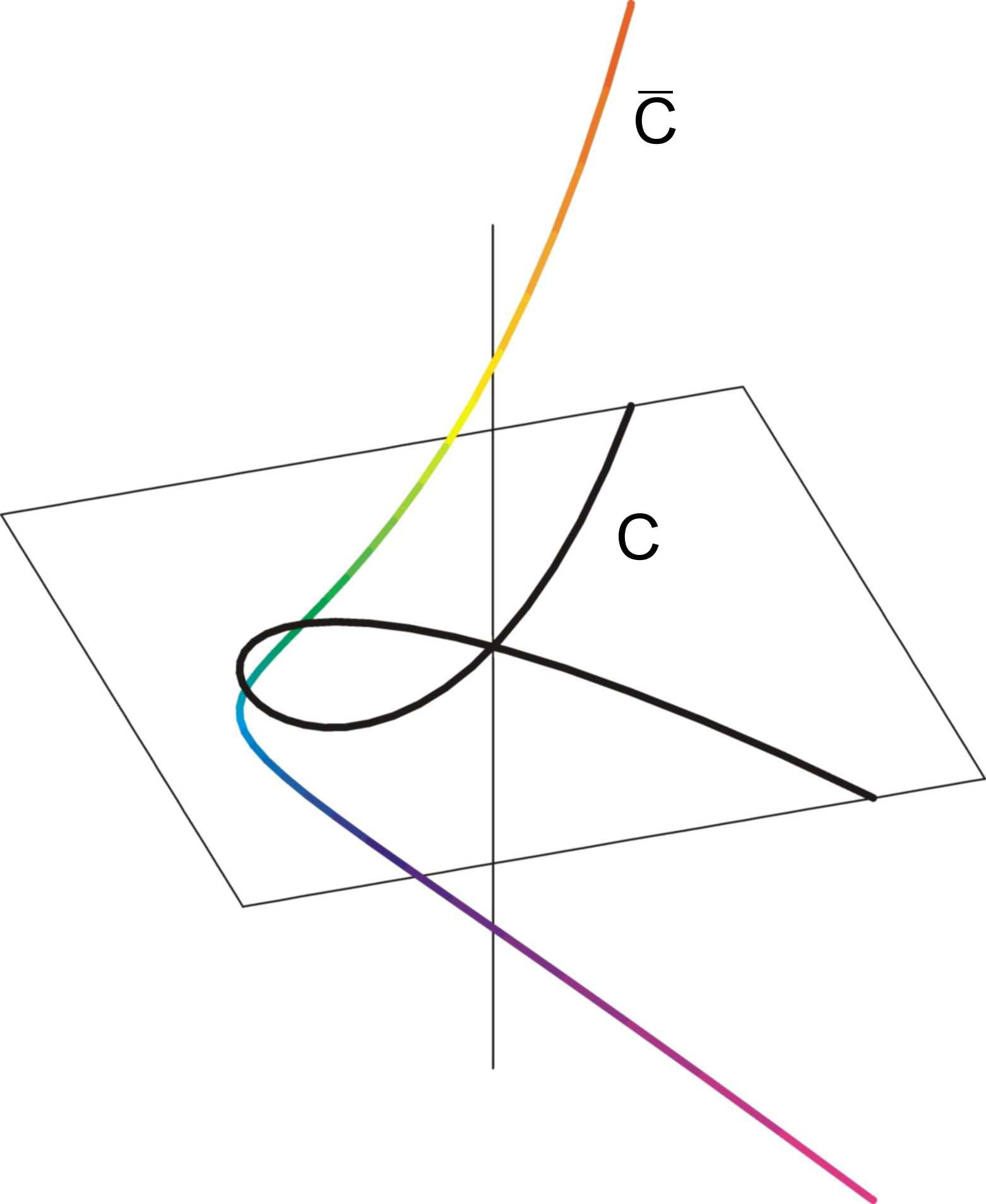}%
\caption{The normalization of the nodal plane curve $C=V(x^3+x^2-y^2)$ is a variant
of the twisted cubic curve $\overline C$ in 3-space.}%
\label{fig normalization}%
\end{center}
\end{figure}
\end{example}

Now, following~\cite{BDLSS}, we describe how the normalization algorithm
can be redesigned so that it becomes parallel in nature.
For simplicity of the presentation, we focus on the case where
$\operatorname{Sing}(A)$ is a finite set. 
This includes the case where $A$ is the coordinate ring of an algebraic curve.

In the example above, the curve under consideration has just one singularity. 
If there is a larger number of singularities, the normalization algorithm as
discussed so far is global in the sense that it ``improves'' all singularities
at the same time. Alternatively, we now aim at ``improving'' the individual
singularities separately, and then put the individual results together. 
In this local-to-global approach, the local computations can be run
in parallel.
We make use of the following result.

\begin{theorem}[\cite{BDLSS}] Suppose that
$
\operatorname{Sing}(A)=\{P_{1},\ldots,P_{r}\}
$
is finite. Then:
\begin{enumerate}
\item  For each $i$,  let
\[
A\,\subseteq\, B_{i}\,\subseteq\,\overline{A}\vspace{-0.15cm}%
\]
be the intermediate ring obtained by applying the normalization algorithm with $P_{i}$ in place of
$J$. Then
\begin{align*}
(B_{i})_{P_{i}} &  =\overline{A_{P_{i}}}, \text{ and}\\
(B_{i})_{Q} &  =A_{Q}\text{ for all }P_{i}\neq Q\in\operatorname{Spec}(A).\vspace{-0.15cm}%
\end{align*}
We call $B_{i}$ the \textbf{minimal local contribution} to $\overline{A}$ at
$P_{i}$.
\item We have
\[
\overline{A}=B_{1}+\ldots+B_{r}\text{.}%
\]
\end{enumerate}
\end{theorem}
This theorem, together with the local version of the Grauert-Remmert criterion,
whose proof is given in~\cite{BDLSS}, yields an algorithm for normalization which 
is often considerably faster than the global algorithm presented earlier, even if 
the local-to-global algorithm is not run in parallel. The reason for this is that
the cost for ``improving''  just one singularity is in many cases much less than 
that for ``improving'' all singularities at the same time. The new algorithm
is implemented in the \textsc{Singular} library \texttt{locnormal.lib}~\cite{locnormal}.
Over the rationals, the algorithm becomes even more powerful by combining it
with a modular approach. This version of the algorithm is implemented in the 
\textsc{Singular} library \texttt{modnormal.lib}~\cite{modnormal}.

\section{Computing GIT-Fans}\label{sec git}

In this section, we give an example of an algorithm that uses Gr\"obner bases, polyhedral computations and algorithmic group theory. It is also suitable for parallel computations.

Recall that one of the goals of Geometric Invariant Theory (GIT)  is to assign to a given algebraic variety $X$ that comes with the action of an algebraic group $G$ in a sensible manner a quotient space $X\good G$.
This setting frequently occurs when we face a variety
$X$ parameterizing a class of geometric objects, for example algebraic curves,
and an action of a group $G$ on $X$ emerging from isomorphisms between 
the objects. There are two main problems.
The first problem is that the homogeneous space $X/G$ is not a good candidate for $X\good G$ as it does not necessarily carry the structure of an algebraic variety.
One then defines for affine $X$ the quotient $X\good G$ as the spectrum of the (finitely generated) invariant ring of the functions of $X$; for general $X$, one glues together the quotients of an affine covering. 
Now a second problem arises: the full quotient $X\good G$ may not carry much information:
For instance,  consider the action of $\mathbb{C}^*:=\CC\setminus \{0\}$ on $X=\mathbb{C}^2$ given by component-wise multiplication
\begin{align}
\mathbb{C}^*\times X\,\to\, X,\qquad
(t,(x,y))\,\mapsto\,(tx, ty).
\label{eq:gitaction}
\end{align}
Then the quotient $X\good \CC^*$ is isomorphic to a point.
However, considering the open subset $U:=X\setminus \{(0,0)\}$ gives us $U\good \CC^* = \PP^1$, the projective line.
For general $X$, there are many choices for these open subsets $U\subseteq X$, where different choices lead to
different quotients 
$U\good G$. To describe this behaviour, Dolgachev and Hu
\cite{DolgachevHu} introduced the \textbf{GIT-fan}, a polyhedral fan describing this variation of GIT-quotients. Recall that a \textbf{polyhedral fan} is a
finite collection of strongly convex rational polyhedral cones such 
that their faces are again elements of the fan and the intersection of any two cones is a common face.

\begin{figure}[H]
\begin{center}
{\includegraphics[width=1.5in]{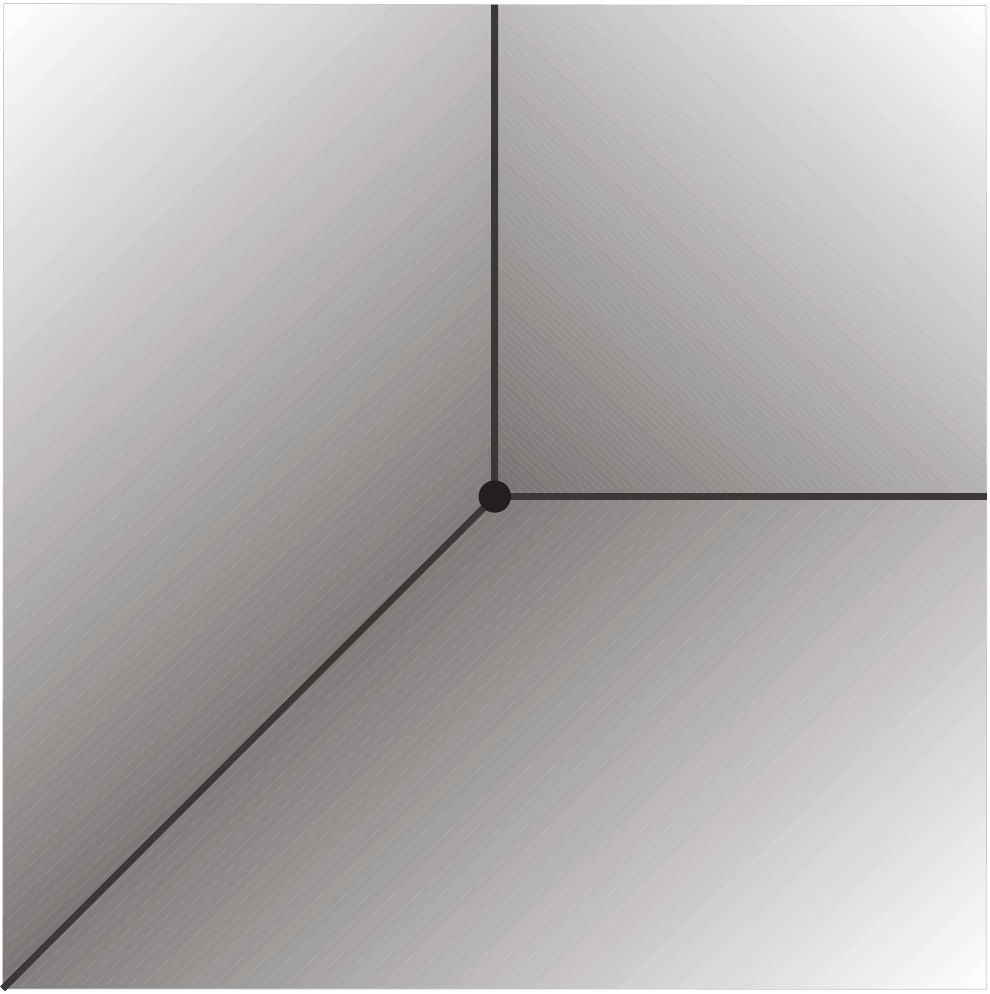}}
\end{center}
\par
\caption{A polyhedral fan in $\mathbb R^2$.}%
\label{fig:fan}
\end{figure}

Of particular importance is the action of an algebraic torus $G=(\mathbb{C}^{\ast})^{k}$,
on an affine variety
$X\subseteq\mathbb{C}^{r}$. In this case, Berchthold/Hausen and
the third author~\cite{BeHa,Ke} have developed a method for computing the GIT-fan, see
Algorithm~\ref{alg:gitfan}. The input of the algorithm consists of

\begin{itemize}[leftmargin=*,itemsep=2pt]

\item an ideal $\mathfrak{a}\subseteq\mathbb{C}[T_{1},\ldots,T_{r}]$ which
defines $X$ and

\item a matrix $Q=(q_{1},\ldots,q_{r})\in\mathbb{Z}^{k\times r}$ such that
$\mathfrak{a}$ is homogeneous with respect to the multigrading defined by setting $\deg(T_i):=q_i\in \mathbb{Z}^k$.
\end{itemize}
Note that the matrix $Q$ encodes the action of  $(\mathbb{C}^{\ast})^{k}$ on $X$. For instance, 
the action~\eqref{eq:gitaction} is encoded in $Q=(1,1)$.

\begin{figure}[h]
\begin{center}
 \includegraphics[width=3.75cm,height=2.5cm]{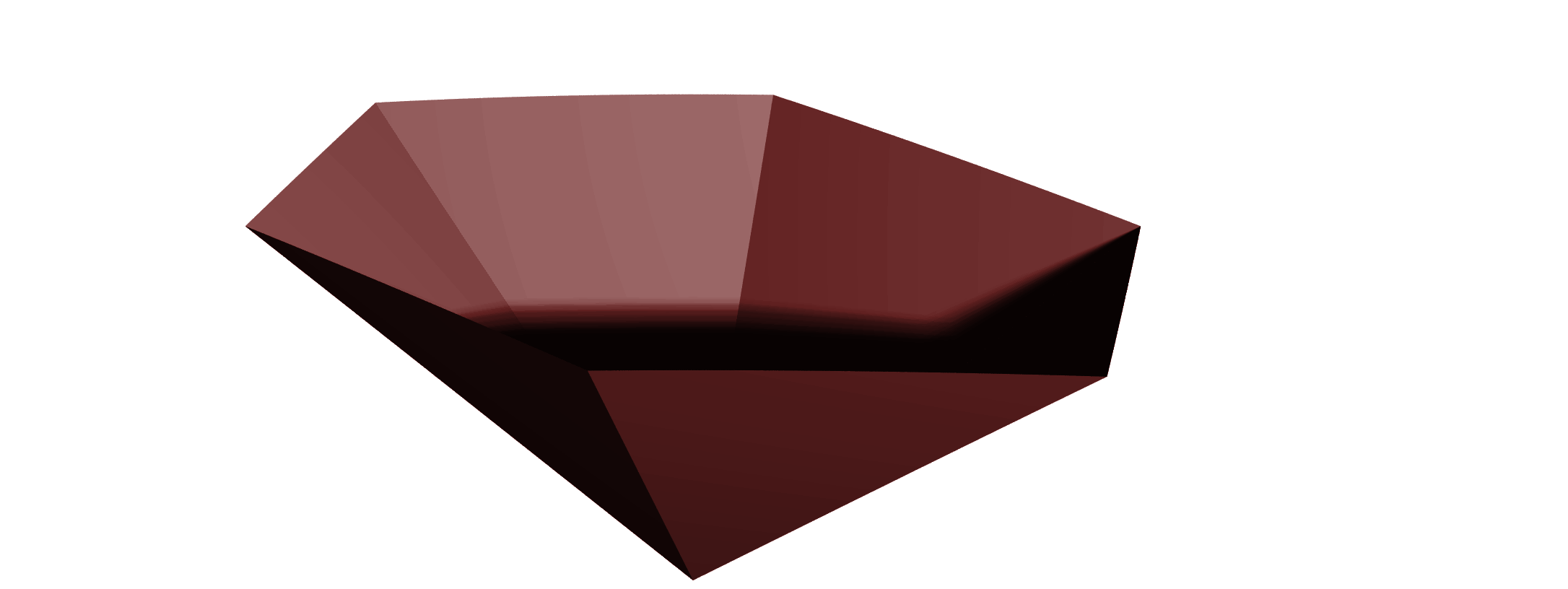}
 \quad
  \includegraphics[width=3.75cm,height=2.5cm]{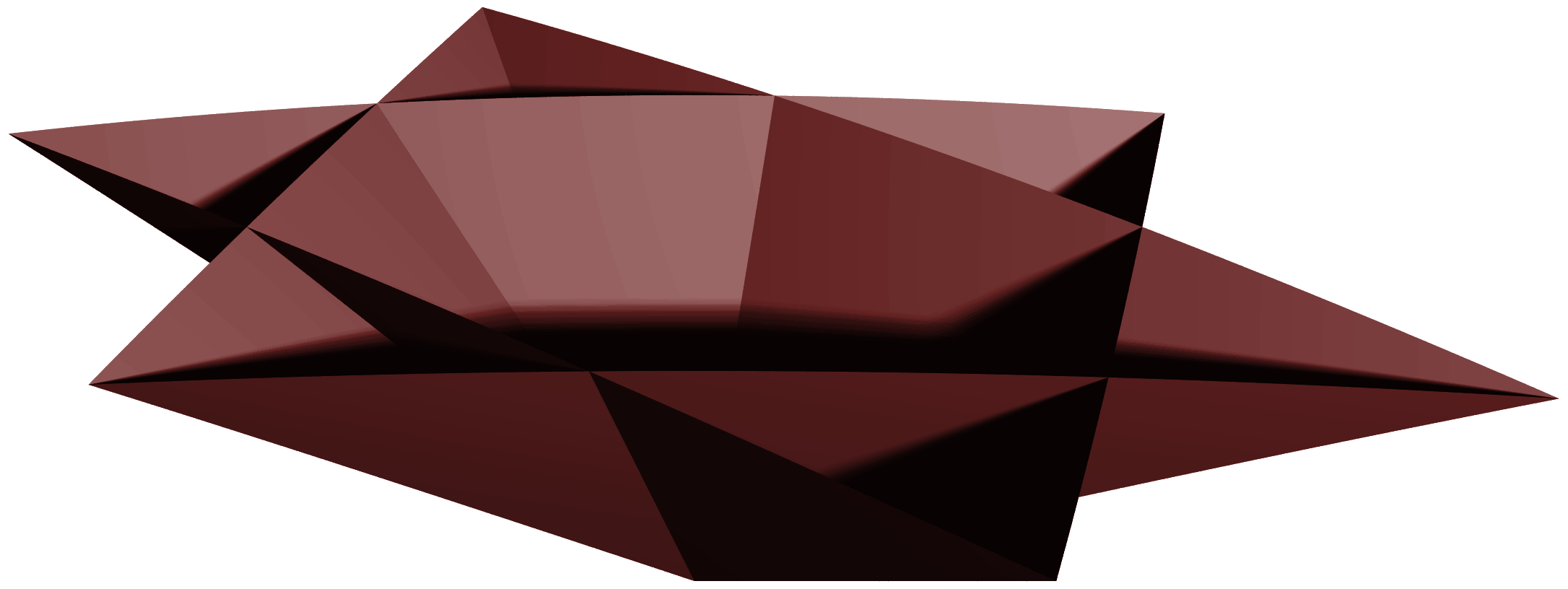}
  \quad
   \includegraphics[width=3.75cm,height=2.5cm]{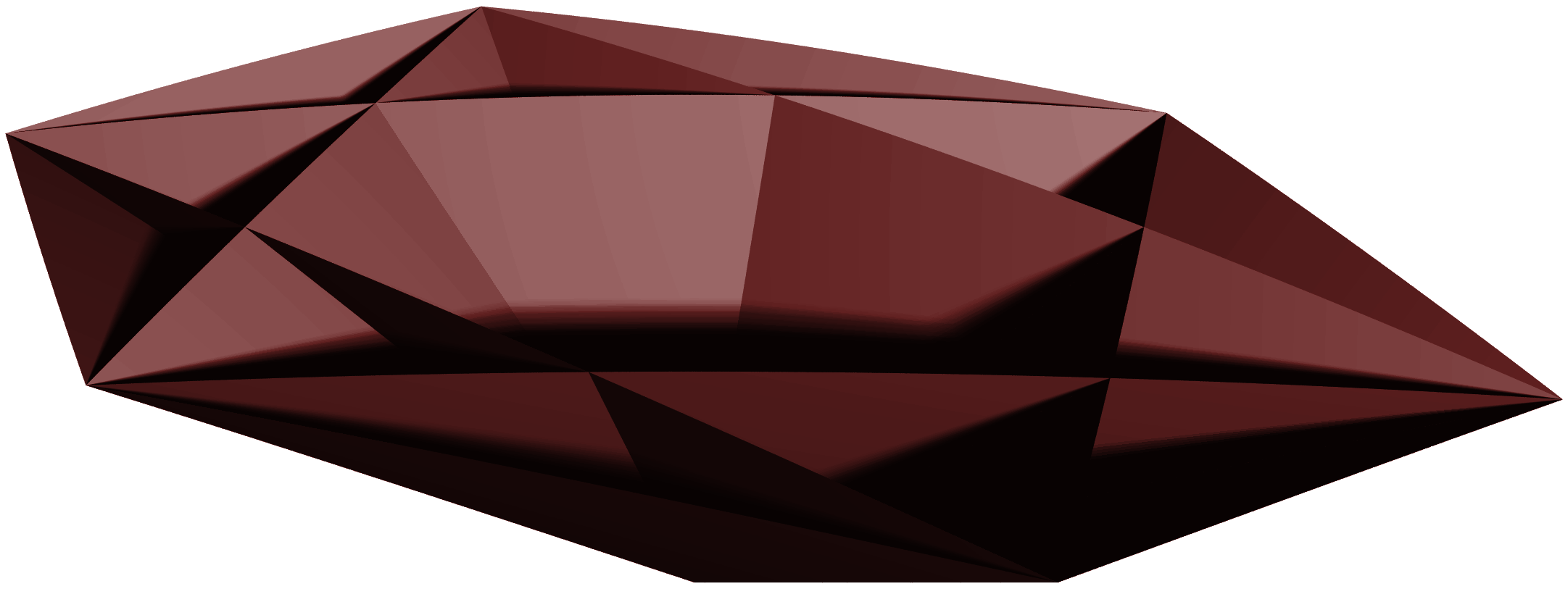}
\end{center}
\caption{Fan traversal.}
\label{fig:traversal}
\end{figure}

Algorithm~\ref{alg:gitfan} can be divided into three main steps.
For the first step, we decompose $\mathbb{C}^{r}$
 into the $2^{r}$ disjoint torus orbits
   \begin{gather*}
    \mathbb{C}^{r}\,=\,\bigcup_{\gamma\subseteq\{1,\ldots,r\}} \!\! O(\gamma),
    \qquad
    O(\gamma)\,:=\,\{(z_{1},\ldots,z_{r})\in\mathbb{C}^{r}\mid z_{i}\neq0 \Leftrightarrow i\in\gamma\}.
   \end{gather*}
The algorithm then identifies in line \ref{afaces} which of the torus orbits $O(\gamma)$ have a
non-trivial intersection with $X$.
The corresponding
$\gamma\subseteq\{1,\ldots,r\}$ (interpreted as faces of the positive orthant $\mathbb{Q}_{\geq 0}^r$)
are referred to as \textbf{$\mathfrak{a}$-faces}. Using the equivalence
\[
X\cap O(\gamma)\,\neq\,\emptyset\quad\Longleftrightarrow\quad\left( \mathfrak{a}%
|_{T_{i}=0 \text{ for } i\notin\gamma}\right)  :\langle T_{1}\cdots
T_{r}\rangle^{\infty}\neq\langle1\rangle,
\]
the $\mathfrak{a}$-faces can be determined by computing the saturation through
Gr\"obner basis techniques available in \textsc{Singular}.
In the second step (line \ref{projectafaces} of the algorithm), the
$\mathfrak{a}$-faces are projected to cones in $\mathbb{Q}^{k}$.
For each $\mathfrak{a}$-face
$\gamma$, defining inequalities and equations of the resulting \textbf{orbit
cones}
\[
Q(\gamma):=\operatorname{cone} (q_{i}\mid i\in\gamma)\,\subseteq\,
\Gamma:=\operatorname{cone}(q_{1},\ldots,q_{r})\, \subseteq\,
\mathbb{Q}^{k}%
\]
are determined, where by $\operatorname{cone}(v_1,\ldots,v_k)$
we mean the polyhedral cone obtained by taking all non-negative linear combinations of the $v_i$.
Computationally, this can be done via the double description method available
in \textsc{polymake}.
We denote by $\Omega$ the set of all orbit cones.
In the final step, the GIT-fan is obtained as
\[
\Lambda(\mathfrak{a}, Q)\,:=\,\{\lambda_{\Omega}(w)\mid w\in\Gamma\} \quad\text{ where }\quad
\lambda_{\Omega}(w) \,:=\, \bigcap_{w\in\eta\in\Omega} \eta.
\]
To compute $\Lambda(\mathfrak{a}, Q)$, we perform a fan-traversal in the following way: Starting with a random
maximal GIT-cone $\lambda_{\Omega}(w_{0})\in \Lambda(\mathfrak{a}, Q)$, we compute its
facets, determine the GIT-cones $\lambda_{\Omega}(w)$ adjacent to it, and
iterate until the support of the fan equals $\operatorname{cone}%
(q_{1},\ldots,q_{r})$. Figure~\ref{fig:traversal} illustrates three steps in such a process.

In line \ref{symdiff} of Algorithm \ref{alg:gitfan}, we
write $\ominus$ for the symmetric difference in the first component.
Again, computation of the facets of a given cone is available through the convex hull
algorithms in \textsc{polymake}. 
\begin{algorithm}[H]
\caption{GIT-fan}
\label{alg:gitfan}
\begin{algorithmic}[1]
\REQUIRE An ideal $\mathfrak{a}\subseteq \mathbb{C}[T_1,\ldots,T_{r}]$ and a matrix $Q\in\mathbb{Z}^{k\times r}$ of full rank such that $\mathfrak{a}$ is homogeneous with respect to the multigrading given by $Q$.
\ENSURE The set of maximal cones of $\Lambda(\mathfrak{a}, Q)$.
\STATE $\mathcal A := \{\gamma\subseteq \mathbb{Q}_{\geq 0}^r \text{ face}\mid \gamma \text{ is an }\mathfrak{a}\text{-face}\}$\label{afaces}
\STATE $\Omega:=\left\{Q (\gamma) \mid \gamma \in  \mathcal A\right\}$\label{projectafaces}
\STATE Choose a vector $w_0\in Q(\gamma)$ such that $\dim(\lambda_\Omega(w_0))=k$.
\STATE $\mathcal C := \{\lambda(w_0)\}$
\STATE $\mathcal{F}:=\{(\tau, \lambda_\Omega(w_0))\mid \tau\subseteq\lambda_\Omega(w_0) \text{ facet with } \tau\not\subseteq\partial\Gamma \}$.
\WHILE{there is $(\eta,\lambda)\in\mathcal{F}$}\label{loop}
\STATE Find $w \in Q(\gamma)$ such that $w\not\in\lambda$ and $\lambda_\Omega(w) \cap \lambda=\eta$.
\STATE $\mathcal C := \mathcal C \cup \{\lambda_\Omega(w)\}$
\STATE $\mathcal F := \mathcal F \ominus \{(\tau,\lambda_\Omega(w))\mid\tau\subseteq\lambda_\Omega(w)  \text{ facet with } \tau\not\subseteq\partial\Gamma\}$\label{symdiff}
\ENDWHILE
\RETURN $\mathcal C$
\end{algorithmic}
\end{algorithm}

Algorithm \ref{alg:gitfan} is implemented in the \textsc{Singular} library
\texttt{gitfan.lib}~\cite{gitfanlib}.
The \textsc{Singular} to
\textsc{polymake} interface \texttt{polymake.so}~\cite{polymakeso}
provides key convex geometry functionality in the \textsc{Singular}
interpreter through a kernel level interface written in C++. We illustrate the
use of this interface by a simple example.

\begin{example}
We compute the normal fan $F$ of the Newton polytope $P$ of the polynomial
$f=x^{3}+y^{3}+1$, see Figure~\ref{fig:fan}. Note that $F$ is the Gr\"obner fan of the ideal $\langle f
\rangle$ and its codimension one skeleton is the tropical 
variety of~$\langle f \rangle$. \medskip

\hspace{-3mm}{ \texttt{$\color{blue}>$ \color{red}LIB \color{black}"polymake.so";} }

\hspace{-3mm}{ \hspace{0.5cm}\texttt{\color{blue}Welcome to polymake version 2.14} }

\hspace{-3mm}{ \hspace{0.5cm}\texttt{\color{blue}Copyright (c) 1997-2015} }

\hspace{-3mm}{ \hspace{0.5cm}\texttt{\color{blue}Ewgenij Gawrilow, Michael Joswig (TU
Berlin)} }

\hspace{-3mm}{ \hspace{0.5cm}\texttt{\color{blue}http://www.polymake.org} }

\hspace{-3mm}{ \hspace{0.5cm}\texttt{\color{blue}// ** loaded polymake.so} }

\hspace{-3mm}{ \texttt{$\color{blue}>$ \color{red}ring \color{black}R = 0,(x,y),dp;
\color{red}poly \color{black}f = x3+y3+1;} }

\hspace{-3mm}{ \texttt{$\color{blue}>$ \color{red}polytope \color{black}P =
newtonPolytope(f); } }

\hspace{-3mm}{ \texttt{$\color{blue}>$ \color{red}fan \color{black} F = normalFan(P);
F;\vspace{0.1cm}} }

\hspace{-3mm}{
\begin{tabular}
[c]{lllll}%
\hspace{0.1cm}%
\begin{tabular}
[c]{l}%
\texttt{\color{blue}RAYS:}\\
\texttt{\color{blue}-1 -1 \#0}\\
\texttt{\color{blue}\ 0\ \ 1 \#1}\\
\texttt{\color{blue}\ 1\ \ 0 \#2}\\
\texttt{\color{blue}MAXIMAL\_CONES:}\\
\texttt{\color{blue}\{0 1\} \#Dimension 2}\\
\texttt{\color{blue}\{0 2\}}\\
\texttt{\color{blue}\{1 2\}}%
\end{tabular}
&  &  &  & \hspace{4cm}\raisebox{-0.567in}{\includegraphics[
height=1.1407in,
width=1.1357in
]{P2fan2-eps-converted-to.pdf}}
\end{tabular}
}
\end{example}

For many relevant examples, the computation of GIT-fans is challenged not only
by the large amount of computations in lines \ref{afaces} and \ref{loop} of Algorithm \ref{alg:gitfan}, but
also by the complexity of each single computation in some boundary cases. 
Making use of symmetries and parallel computations, we can open up the
possibility to handle many interesting new cases by considerably simplifying
and speeding up the computations. 
For instance, the computations in line \ref{afaces} of
Algorithm \ref{alg:gitfan} can be executed independently in parallel. Parallel
computation techniques can also be applied in the computation of
$\lambda_{\Omega}(w_{0})$ and the traversal of the GIT-fan. This step, however, is not trivially parallel.

An example of the use of symmetries is~\cite{BKR}; here, the first, third and fourth authors have applied and extended the technique described above to obtain the cones of the Mori chamber decomposition (the GIT-fan of the action of the characteristic torus on its total coordinate space) %
of the Deligne-Mumford compactification $\overline
{M}_{0,6}$ of the moduli space of $6$-pointed stable curves of genus zero that lie within the cone of movable divisor classes. A
priori, this requires to consider $2^{40}$ torus orbits in line~\ref{afaces}.
Hence, a direct application of Algorithm~\ref{alg:gitfan} in its stated form
is not feasible. However, moduli spaces in algebraic geometry often have large
amounts of symmetry. For example, on $\overline{M}_{0,6}$ there is a natural
group action of the symmetric group~$S_{6}$
which Bernal~\cite{Ber} has extended to the input data $\mathfrak{a}$ and $Q$ required for Algorithm~\ref{alg:gitfan}.
The GIT-fan $\Lambda(\mathfrak{a}, Q)$,
and all data that arises in its computation reflect these symmetries. 
Hence, by computing an orbit decomposition under the action of the group
of symmetries of the set of all torus orbits, we
can restrict ourselves to a distinct set of representatives. Also the 
fan-traversal can be done modulo symmetry.
 To compute the orbit decomposition, we
apply the algorithms for symmetric groups implemented in \textsc{GAP}.

\begin{example}
We apply this technique in the case of the affine cone $X$ over the Grassmannian $\mathbb{G}(2,5)$ of
$2$-dimensional linear subspaces in a $5$-dimensional vector space, see also~\cite{BKR}. By making
use of the action of $S_{5}$, the number of monomial containment tests in line
\ref{afaces} can be reduced from $2^{10}=1024$ to $34$. A distinct set of
representatives of the orbits of the $172$ $\mathfrak{a}$-faces consists 
of $14$ elements. The GIT-fan has $76$ maximal cones, which fall into $6$
orbits. Figure \ref{fig G25} shows both the adjacency graph of the maximal cones of
the GIT-fan and that of their orbits under the $S_{5}$-action.
This GIT-fan has also been discussed in~\cite{Ber,DolgachevHu,AH}. Note that by 
considering orbits of cones not only the computation of the fan is considerably simplified,
but also the theoretical understanding of the geometry becomes  easier.
\begin{figure}[h]
\begin{center}
\includegraphics[width = 0.8\textwidth]{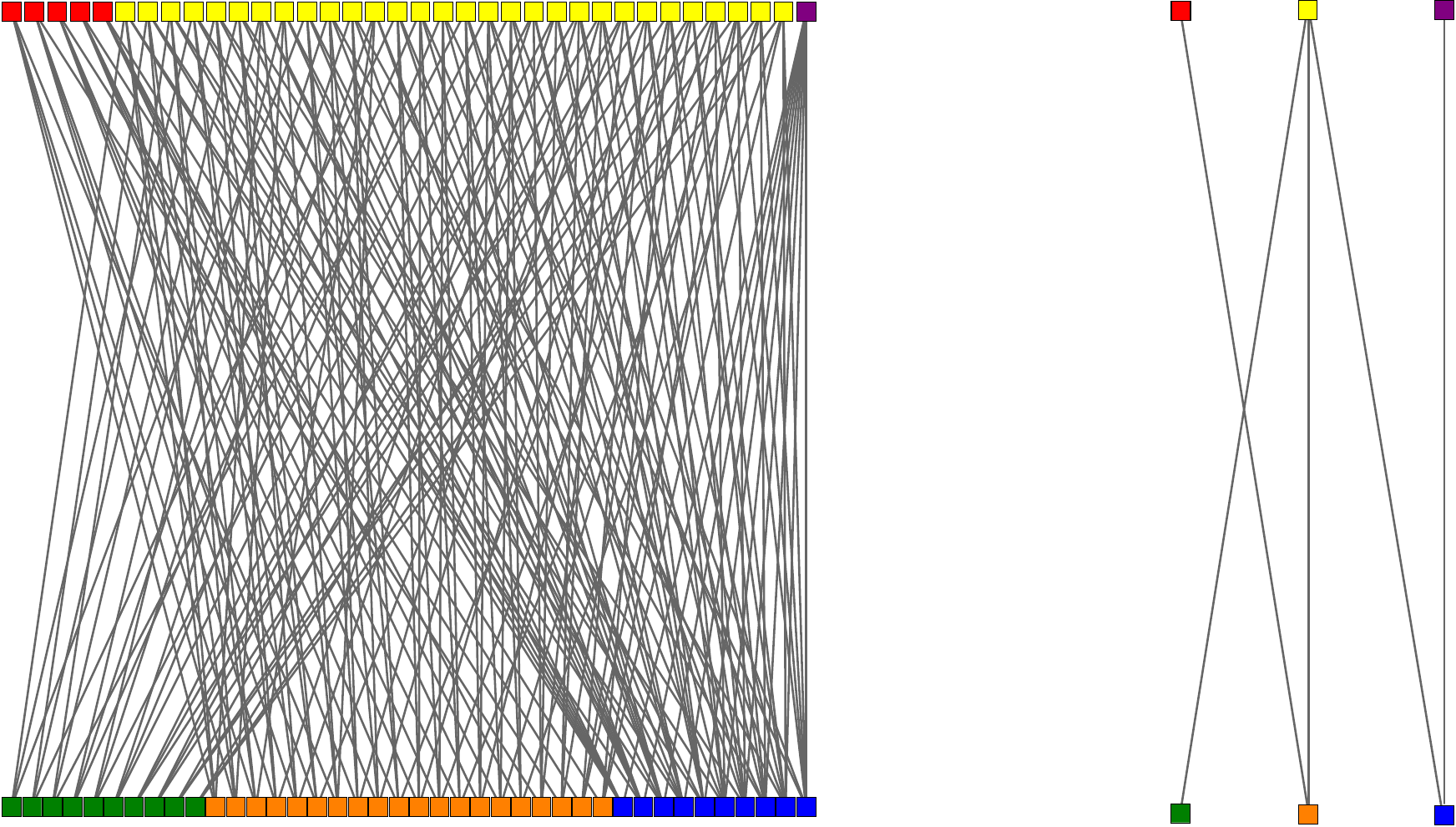}
\end{center}
\par
\caption{The adjacency graph of the set of maximal cones of the
GIT-fan of $\mathbb{G}(2,5)$ and the adjacency graph of the orbits of these cones
under the $S_{5}$-action.}%
\label{fig G25}
\end{figure}
\end{example}

To summarize, Algorithm \ref{alg:gitfan} requires the following key
computational techniques from commutative algebra, convex geometry, and group theory:

\begin{itemize}[leftmargin=*,itemsep=2pt]

\item Gr\"obner basis computations,

\item convex hull computations, and

\item orbit decomposition.
\end{itemize}

These techniques are provided by 
\textsc{Singular}, \textsc{polymake}, and \textsc{GAP}. At the current stage, \textsc{polymake} can be used from \textsc{Singular} in a convenient way through \texttt{polymake.so}. An interface to use \textsc{GAP} functionality directly from \textsc{Singular} is subject to future development.

\bibliographystyle{abbrv}
\bibliography{bdkr}

\end{document}